\documentclass{amsart}
\usepackage{amssymb}
\usepackage{mathrsfs}
\usepackage{stmaryrd} 
\usepackage{bbm}
\usepackage{oldgerm}
\usepackage[francais]{babel}
\usepackage[all]{xy}

\usepackage{verbatim, amsmath, amscd, amssymb}

\newtheorem{teo}[subsection]{Th\'{e}or\`{e}me}
\newtheorem{prop}[subsection]{Proposition}
\newtheorem{cor}[subsection]{Corollaire}
\newtheorem{lem}[subsection]{Lemme}

\theoremstyle{definition}

\newtheorem{defi}[subsection]{D\'{e}finition}
\newtheorem{rema}[subsection]{Remarque}

\numberwithin{equation}{subsection}

\newcommand{\mN}{{\mathbb N}}

\newcommand{\mZ}{{\mathbb Z}}

\newcommand{\cA}{{\mathcal A}}
\newcommand{\cB}{{\mathcal B}}
\newcommand{\cC}{{\mathcal C}}

\newcommand{\cL}{{\mathcal L}}

\newcommand{\cO}{{\mathcal O}}

\newcommand{\bbF}{\mathbb F}

\newcommand{\bbN}{\mathbb N}
\newcommand{\bbQ}{\mathbb Q}

\newcommand{\bbZ}{\mathbb Z}

\newcommand{\bfU}{\mathbf U}
\newcommand{\bU}{{\bf U}}

\newcommand{\Fr}{\mathrm{Fr}}
\newcommand{\Frac}{\mathrm{Frac}}

\newcommand{\rH}{\mathrm{H}}

\newcommand{\rR}{\mathrm{R}}

\newcommand\pf{\noindent {\bf Proof:  }}

\newcommand{\id}{\mathrm{id}}

\newcommand{\ind}{\mathrm{ind}}
\newcommand{\lbr}{\begin{bmatrix}}
\newcommand{\rbr}{\end{bmatrix}}
\newcommand{\St}{\mathrm{St}}
\newcommand{\Dist}{\mathrm{Dist}}

\newcommand{\Mod}{\mathbf{Mod}}

\begin{document}

\title[Un scindage du morphisme de Frobenius quantique]
{Un scindage du morphisme de Frobenius quantique*}
 \author{ Michel Gros et Masaharu Kaneda}
 
\address{M.G. CNRS UMR 6625, IRMAR, Universit\'{e} de Rennes 1,
Campus de Beaulieu, 35042 Rennes cedex, France}
\address{M.K. Osaka City University, Department of Mathematics,
3-3-138 Sugimoto, Sumiyoshi-ku, Osaka 558-8585, Japan}
\begin{abstract}
Nous montrons que le morphisme de Frobenius quantique construit par Lusztig dans le cadre des alg\`ebres enveloppantes quantiques $U_\cB$ sp\'ecialis\'ees en une racine de l'unit\'e   admet un scindage multiplicatif (non unitaire). Nous utilisons pour ce faire  une base de la partie torique de la petite alg\`ebre quantique constitu\'ee d'idempotents orthogonaux deux \`a deux et de somme 1 et faisons de m\^eme dans le cas   ``modulaire'' pour l'alg\`ebre des distributions d'un groupe alg\'ebrique semi-simple.  

We show that the quantum Frobenius morphism constructed by Lusztig in the setting of the quantum enveloping algebra $U_\cB$ specialized at a root of unity  admits a multiplicative splitting (non unital). We construct the splitting using a basis of the toral part of 
the small quantum algebra consisting of pairwise orthogonal idempotents summing up to 1, and likewise in the modular case
of the algebra of distributions for a semisimple algebraic group.
\end{abstract}
  
\email{michel.gros@univ-rennes1.fr}
\email{kaneda@sci.osaka-cu.ac.jp}

\thanks{* supported in part by JSPS Grants in Aid for Scientific Research 23540023}

\maketitle

\setcounter{tocdepth}{1}
\tableofcontents

\section{Introduction}
\subsection{}\label{intro1}

Soient $A$ une matrice de Cartan  $\ell\times \ell$ de type fini, 
$\cA=\bbZ[v,v^{-1}]$ l'anneau des polyn\^omes de   Laurent  en l'ind\'etermin\'ee   $v$, et $U$ la $\cA$-forme (i.e. avec ``puissances divis\'ees'') de l'alg\`ebre quantique associ\'ee \`a $A$.
Soit maintenant $l>1$ un entier   premier avec tous les coefficients de
$A$.
Soient $\cB$ le quotient de $\bbZ[\frac{1}{2}][v]$ par l'id\'eal engendr\'e par  le $l$-i\`eme polyn\^ome  cyclotomique $\Phi_{l}$, $q$ l'image de $v$ dans $\cB$, et $U_\cB=U\otimes_\cA\cB$.
Soient enfin $\bU$ la $\bbZ$-forme de Kostant de l'alg\`ebre enveloppante universelle associ\'ee \`a $A$,
et
$\bU_\cB=\bU\otimes_\bbZ
\cB$.
  Lusztig a \'etabli  dans \cite[8.10]{lus}    l'existence d'un
morphisme de Frobenius quantique  
$\Fr: U_\cB\to \bU_\cB$,
qui, lorsque $l$ est premier et not\'e alors $p$,  est un rel\`evement du morphisme de    Frobenius usuel sur l'alg\`ebre des 
distributions correspondante sur le corps premier $\bbF_{p}$ \`a    $p$ \'el\'ements.
Nous construisons dans cet article un homomorphisme {\emph{non unif\`ere}} d'alg\`ebres (\cite[III. \S 1, 1]{Bour}) $\phi:\bU_\cB \to U_\cB$  scindant  (cf.  \ref{teointrod})   le morphisme    $\Fr$   et relevant,    en un sens qu'on pr\'ecisera  \eqref{relev}, le scindage du Frobenius usuel construit dans \cite[Thm. 1.3]{gk}.

\subsection{}\label{intro2}
Dans une situation similaire mais avec toutefois   des hypoth\`eses   diff\'erentes sur $l$ (cf. \cite[\S 2 ]{Li}),  adaptant donc en cons\'equence    les notations pr\'ec\'edentes,   rappelons  que 
Littelmann a   d\'efini  \`a l'aide de
\cite[Thm 1]{Li}
la  {\em{contraction}} d'un $U_\cB$-module de Weyl  
 dans le but de compl\'eter son programme visant \`a \'etablir une th\'eorie 
de  ``mon\^omes standards'' pour  les $\bU_\cB$-modules de Weyl duaux.
 Il utilise simplement pour ce faire  un certain scindage sur la seule partie nilpotente de  $\bU_\cB$.  Le prolongement  na\"if de son scindage \`a tout 
$\bU_\cB$ 
s'av\`ere ne  pas \^etre multiplicatif. 
N\'eanmmoins,  revenant \`a  nos hypoth\`eses \ref{intro1}, reprenant sa construction et   faisant intervenir  en plus
une mesure involutive invariante de
la partie torique de la sous-alg\`ebre infinit\'esimale de $U_\cB$, un miracle analogue \`a 
celui du cas modulaire  \cite{gk} nous permet de construire  un prolongement {\emph{multiplicatif}}.
 On est alors en mesure de contracter n'importe quel 
 $U_\cB$-module pour en faire un 
  $\bU_\cB$-module.
Si bien s\^ur on ne s'int\'eresse qu'aux $U_\cB$-modules ayant une d\'ecomposition suivant leurs poids,
on   dispose  \cite[Prop. 3.4]{Mc} d\'ej\`a d'une fa\c{c}on de  les contracter : il suffit d'interpr\'eter ces modules    comme des
modules unitaires (``unital modules'') sur l'alg\`ebre quantique modifi\'ee pour laquelle un
scindage du Frobenius est d\'efini dans \emph{loc. cit.} (et ce, sans m\^eme la restriction sur $l$ ci-dessus). 
 
\subsection{}\label{intro3}
Pour \'enoncer plus pr\'ecis\'ement notre r\'esultat principal,
notons $\bbQ(v)$ (resp. $\bbQ(q)$)
le corps de fractions de  $\cA$ (resp. $\cB$) et
 $U_{\bbQ(v)}$ l'alg\`ebre quantique associ\'ee \`a  $A=[\!(a_{ij})\!]$ sur
$\bbQ(v)$  avec  ses g\'en\'erateurs standards
$E_i, K_i, F_i$, $i\in[1,\ell]$.
Soient $d_i\in\{1,2,3\}$, $i\in[1,\ell]$, tels que la matrice
$[\!(d_ia_{ij})\!]$ soit sym\'etrique, et
$v_i=v^{d_i}$.
Pour tous $i\in[1,\ell]$ et $n\in\bbN$
posons
$[n]_i^!=\frac{v_i^n-v_i^{n}}{v_i-v_i^{-1}}$,  
$E_i^{(n)}=\frac{E_i^n}{[n]_i^!}$,
$F_i^{(n)}=\frac{F_i^n}{[n]_i^!}$.
Alors, $U$ est simplement la $\cA$-sous-alg\`ebre de   $U_{\bbQ(v)}$
engendr\'ee par les 
$E_i^{(n)}, F_i^{(n)}, K_i, K_i^{-1}$,
$i\in[1,\ell]$,
$n\in\bbN$.
Soient \'egalement
$X_i, Y_i$, $i\in[1,\ell]$, 
les g\'en\'erateurs standards de l'alg\`ebre enveloppante 
universelle 
$\bU_{\bbQ(q)}$
associ\'ee \`a $A$ sur
$\bbQ(q)$.
Si
 $X_i^{(n)}=\frac{X_i^n}{n!}$
et
$Y_i^{(n)}=\frac{Y_i^n}{n!}$,
$i\in[1,\ell]$,
$n\in\bbN$, 
$\bU_\cB$ est la  $
\cB$-sous-alg\`ebre de
$\bU_{\bbQ(q)}$ 
engendr\'ee par les
$X_i^{(n)}$ et
$Y_i^{(n)}$,
pour $i\in[1,\ell]$ et
$n\in\bbN$.
Le morphisme de Frobenius introduit par Lusztig  (strictement dit, Lusztig travaille  dans \cite{lus} au-dessus de 
$\bbZ[v]$ modulo l'id\'eal engendr\'e par $\Phi_{l}$)
\begin{equation}\label{Fr} 
\Fr:U_\cB\to\bU_\cB
\end{equation} 
est alors d\'efini,   pour tous  $i\in[1,\ell]$,  $n\in\bbN$, par
\begin{equation}\label{defFr} 
E_i^{(n)}\mapsto
\begin{cases}
X_i^{(\frac{n}{l})}
&\text{si $l|n$}
\\
0
&\text{sinon}  
\end{cases}
\
,
\quad
F_i^{(n)}\mapsto
\begin{cases}
Y_i^{(\frac{n}{l})}
&\text{si $l|n$} 
\\
0
&\text{sinon}  
\end{cases},
\quad
K_i^{\pm1}\mapsto1.
\quad
\end{equation}
\subsection{}\label{intro4}
Si $\bU_\cB^{\pm}$ d\'esigne la  $
\cB$-sous-alg\`ebre de 
$\bU_\cB$ engendr\'ee par 
les
$X_i^{(n)}$
(resp.
$Y_i^{(n)}$),
pour $i\in[1,\ell]$ et $ n\in\bbN$,
Lusztig
 d\'efinit  \cite[8.6]{lus} un scindage 
$\Fr'^{\pm}$ de
$\Fr\mid_{U_\cB^{\pm}}$
par
$X_i^{(n)}\mapsto
E_i^{(nl)}$
(resp.
$Y_i^{(n)}\mapsto
F_i^{(nl)}$)
pour tous
$i\in[1,\ell]$
et 
$n\in\bbN$.
Si
$\bU_\cB^{\geq0}$
(resp.
$\bU_\cB^{\leq0}$)
d\'esigne la $
\cB$-sous-alg\`ebre de
 $\bU_{\bbQ(q)}$ engendr\'ee par
$\bU_\cB^{\pm}$
et par les $\binom{H_i}{n}=\frac{H_i(H_i-1)\dots(H_i-n+1)}{n!}$,
$i\in[1,\ell]$,
$n\in\bbN$, 
Kumar et Littelmann 
\cite[Thm. 1.2]{kl}
prolongent ensuite 
$\Fr'^{\pm}$ \`a
$\bU_\cB^{\geq0}$ et 
$\bU_\cB^{\leq0}$, mais  au prix du passage au quotient (lequel est indispensable
 pour conserver la multiplicativit\'e) par l'id\'eal
$(K_i^l-1\mid
i\in[1,\ell])$.
Si l'on introduit maintenant pour $i\in[1,\ell]$ 
\begin{equation}\label{introkappa}
\kappa_{i0}=\frac{1}{2l} \sum_{j=0}^{2l-1}K_i^{j}\in
U_{\bbQ(v)} \,\,\, ; \,\,\,  \kappa=\kappa_{10}\dots\kappa_{\ell0},
\end{equation} 
 il s'av\`ere que
$\kappa$ appartient en fait \`a 
$U_\cB$ et est l'unique idempotent  non trivial
de la  $\cB$-sous-alg\`ebre 
$u_\cB^0$
de
$U_\cB$ engendr\'ee par les
$K_i$, $i\in[1,\ell]$, \`a \^etre invariant sous l'action 
naturelle \`a gauche, autrement dit :
$\kappa\in\{x\in u_\cB^0\mid
yx=
\varepsilon(y)x\  \,{\text{pour tout}}\, y\in u_\cB^0\}$
(avec 
$\varepsilon$ la co\"unit\'e de
$u_\cB^0$).
Cet \'el\'ement  $\kappa$ apparait donc comme la variante quantique de
la mesure invariante $\mu_0$, introduite dans \cite[1.3]{gk}, 
sur le noyau de Frobenius d'un tore maximal du groupe alg\'ebrique semi-simple $G$ sur ${\Bbb{F}}_{p}$ correspondant. On peut alors \'enoncer le
\begin{teo}\label{teointrod} 
\begin{itemize}
\item[{\rm (i)}]
 Il existe un   unique  homomorphisme (non unif\`ere) de $\cB$-alg\`ebres   \begin{equation}\label{defscind}
\phi:\bU_\cB\to
U_\cB
\end{equation}
 tel  que 
$X_i^{(n)}\mapsto
E_i^{(nl)}\kappa$,
$Y_i^{(n)}\mapsto
F_i^{(nl)}\kappa$
pour tous 
$i\in[1,\ell]$ et $n\in\bbN$. 
\item[{\rm (ii)}]
L'application $\phi$ se factorise \`a travers le facteur direct $\kappa
U_\cB\kappa$ de $U_\cB$  en un homomorphisme unif\`ere de $\cB$-alg\`ebres. On a   
$\Fr\circ\phi=\id_{\bU_\cB}$.

\end{itemize}
\end{teo} 
La d\'{e}monstration  dont nous disposons  est tout \`a fait terre-\`a-terre.   Elle nous a sembl\'e n\'eanmmoins  m\'erit\'ee d'\^etre donn\'ee, le r\'esultat sugg\'erant par exemple une alternative \`a l'emploi des alg\`ebres quantiques modif\'ees pour certaines applications.  L'article est structur\'e en trois grandes parties. Nous traitons tout d'abord en d\'etail (\S2 et \S3)  le cas de l'alg\`ebre $\mathfrak{sl}_2$ en montrant l'int\'egralit\'e de $\kappa=\kappa_{10}$ ainsi que de quelques autres idempotents utiles pour la suite et en \'etablissant l'existence du scindage $\phi$. Nous en d\'eduisons l'existence de   d\'ecompositions associ\'ees de $U_\cB$. La seconde partie (section 4) est d'un int\'er\^et assez largement ind\'ependant. Les idempotents qui sont naturellement apparus lors la d\'emonstration de \ref{teointrod}  invitent \`a  revenir sur leurs analogues modulaires et leurs g\'en\'eralisations (lesquelles sont sp\'ecifiques à la caract\'eristique $p>0$). Nous montrons comment ils permettent de d\'ecrire pr\'ecis\'ement (cf. Prop. \ref{descrquant}, Prop. \ref{descrmodul}) plusieurs des alg\`ebres apparaissant dans ce travail et aussi, comment dans le cas modulaire on peut obtenir une nouvelle  preuve de l'existence (cf. Thm \ref{descrDist})  d'un scindage du morphisme de Frobenius obtenu pr\'ec\'edemment dans \cite{g}.  Enfin,  dans la derni\`ere partie (\S5 et \S6), nous expliquons  les am\'enagements, essentiellement  typographiques,  \`a apporter   pour traiter le cas g\'en\'eral du th\'eor\`eme \ref{teointrod} et terminons par  quelques corollaires laissant d'autres applications pour un travail ult\'erieur. 

Le premier auteur (M.G.) remercie tr\`{e}s sinc\`{er}ement le second (M.K.)  pour son g\'{e}n\'{e}reux accueil  lors de son s\'{e}jour en mai 2012 \`{a} l'Universit\'{e} de la Ville d'Osaka (OCU) ayant rendu possible ce travail   ainsi que le d\'{e}partement de math\'{e}matiques pour l'atmosph\`{e}re  stimulante dans laquelle il a    \'{e}t\'{e} effectu\'{e}. Apr\`es une premi\`ere r\'edaction de ce travail, nous avons pu b\'en\'eficier d'un commentaire  lumineux de Yoshihisa Saito sur celle-ci nous conduisant  \`a    repenser et à simplifier consid\'erablement certaines d\'emonstrations des \S2 et  \S3 : nous lui adressons nos remerciements les plus chaleureux. Nous remercions \'egalement  H.H. Andersen d'avoir attir\'e notre attention sur \cite{PS}.

\section{Le cas  $\mathfrak{sl}_2$ : sur l'int\'egralit\'e  de  quelques idempotents}\label{1}

\subsection{}\label{integ}

On se place ici dans le cas $\ell=1$, $A=(2)$ et $l$ est par cons\'equent un nombre  impair quelconque. Rappelons qu'on a (avec un all\`egement de notations \'evident)  dans  $
U_{\cB}$ les relations
\begin{equation} \label{relgen}
KE = q^{2}EK \,\, ; \,\, KF= q^{-2}FK  \,\, ; \,\, [E,F] = \frac{K-K^{-1}}{q-q^{-1}} .
\end{equation}  
et qu'on note, pour $c\in \mZ$, $m\in\mN$
\begin{equation} \label{binomK}
 \left[\begin{array}{c}K ;c \\m\end{array}\right] = \prod_{h=1}^{m} \frac{Kv^{c-h+1} - K^{-1}v^{-(c-h+1)}}{v^{h}-v^{-h}} \in  
 U_{\bbQ(v)}. 
\end{equation} 
Cet \'el\'ement appartient en fait \`a  
$U$ et nous noterons par le m\^eme symbole son image dans  $U_\cB$. 
Nous poserons aussi simplement
\begin{equation}\label{defbinomK}
 \left[\begin{array}{c}K  \\m\end{array}\right] = \left[\begin{array}{c}K ;0 \\m\end{array}\right]  \in U_\cB. 
\end{equation}

Notons
\begin{equation}\label{defUq}
U_q=U_\cB\otimes_\cB \bbQ(q) \,\,\,; \,\,\,u_{q}^0= u_\cB^{0}\otimes_\cB \bbQ(q)  = \bbQ(q) [K].
\end{equation}

Pour toute la suite de ce travail, nous noterons     $q^{\frac{1}{2}} = -q^{-\frac{l-1}{2}}\in \cB$ qui est une une racine primitive $2l$-i\`eme de 1.

Posons, pour $n\in\bbZ$,  
 \begin{equation}\label{defkappa'j}
 \kappa_n =\frac{1}{2l}\sum_{i=0}^{2l-1}q^{\frac{-ni}{2}}K^i\in
U_q 
 \end{equation}
 On a de mani\`ere \'evidente    $\kappa_0 =\kappa$ \eqref{introkappa}  et $\kappa_n =\kappa_m $ si   $n\equiv m\mod2l$.
 
 On v\'erifie imm\'ediatement \`a l'aide de la relation $K^{2l}=1$  que 
\begin{equation}\label{weightofkappa}
K\kappa _n=q^{{\frac{n}{2}}}\kappa _n.
\end{equation}
Il s'ensuit que les 
$\kappa _n$ pour $n\in[0,2l[$, appartenant \`a des sous-espaces propres deux-\`a-deux distincts pour l'action de $K$, forment une base du   $\bbQ(q)$-espace vectoriel,   
$u_q^0$.
De plus, comme les  $\bbQ(q)\kappa_n$ pour $n\in[0,2l[$ sont des id\'eaux deux-\`a-deux distincts de  $u_q^0$,
les  $\kappa_n$ doivent \^etre orthogonaux entre eux lorsqu'ils sont distincts :  
\begin{equation}\label{}
\kappa _i\kappa _j=0
\end{equation}
pour tout $i\ne j$.

Enfin, utilisant de nouveau que $K^{2l}=1$, on v\'erifie que chaque 
$\kappa _n$ est un  idempotent et que l'on a donc 
\begin{equation}
1=\sum_{n=0}^{2l-1}\kappa _n.
\end{equation}
 \begin{lem}\label{intkappa}
 On a 
\begin{equation}\label{}
  \kappa _{0} = \frac{1}{2}. \{  \sum_{i=0}^{l-1}    (-1)^{i}     \left[\begin{array}{c}K \\i\end{array}\right] (q^{i} +   q^{ - i} K)  \} \in u^0_\cB.
\end{equation}\label{teo1eq2}

 \end{lem}

 Rappelons 
(cf. par exemple \cite{lus2}, preuve de Prop. 2.14 ) d'abord qu'on a, pour tout $i\geq0$,  la relation
\begin{equation}\label{relceccu}
  K^{2}\left[\begin{array}{c}K   \\i\end{array}\right]  =   q^{i}(q^{i+1}-q^{-i-1})   K\left[\begin{array}{c}K    \\i+1\end{array}\right] +   q^{2i}  \left[\begin{array}{c}K \\i\end{array}\right].
  \end{equation}
On en d\'eduit que  
\begin{align}
K (\frac{1}{2}. 
&
\{  \sum_{i=0}^{l-1}    (-1)^{i}    
\lbr
K \\i\rbr (q^{i} +   q^{ - i} K)  \})
=
\frac{1}{2}. \{  \sum_{i=0}^{l-1}    (-1)^{i}     \lbr
K \\i\rbr (q^{i} K+   q^{ - i} K^2)  \}
\\
\notag&=
\frac{1}{2} \sum_{i=0}^{l-1}    (-1)^{i} 
\{    \lbr
K \\i\rbr q^{i} K+   q^{ - i}
(q^i(q^{i+1}-q^{-i-1})K
\lbr
K
\\
i+1
\rbr
+q^{2i} 
\lbr
K
\\
i
\rbr)  \}
\quad
\text{par (\ref{relceccu})}  
\\
\notag&=
\frac{1}{2} \sum_{i=0}^{l-1}    (-1)^{i} 
( \lbr
K \\i\rbr q^{i} K+   (q^{i+1}-q^{-i-1})K
\lbr
K
\\
i+1
\rbr
+q^{i} 
\lbr
K
\\
i
\rbr)  
\\
\notag&=
\frac{1}{2}\{ \sum_{i=0}^{l-1}    (-1)^{i}     \lbr
K \\i\rbr (q^{i} +   q^{ - i} K)  \}+
(q^l-q^{-l})K\lbr
K
\\
l
\rbr
\\
\notag&=
\frac{1}{2}\sum_{i=0}^{l-1}    (-1)^{i}     \lbr
K \\i\rbr (q^{i} +   q^{ - i} K).
\end{align}
Gr\^ace \`a (\ref{weightofkappa}), on obtient donc que 
$\frac{1}{2}\sum_{i=0}^{l-1}    (-1)^{i}     \left[\begin{array}{c}K \\i\end{array}\right] (q^{i} +   q^{ - i} K)\in
\bbQ(q)^\times\kappa_0$.
Posons maintenant
$\kappa=\kappa_0$
et
$\kappa'=
\frac{1}{2}\sum_{i=0}^{l-1}    (-1)^{i}     \left[\begin{array}{c}K \\i\end{array}\right] (q^{i} +   q^{ - i} K)$.
Pour voir que $\kappa' =\kappa$, 
il suffit, puisque $\kappa'\in\bbQ(q)\kappa$,
de voir que
$\kappa'\kappa=\kappa$. 
Si l'on \'ecrit
$\kappa'=\lambda\kappa$, on a
$\kappa'\kappa=\lambda\kappa$.
Mais l'on sait que
$K\kappa=
\kappa$.
Comme (\ref{binomK})
$\lbr
K\\
i\rbr
=
\prod_{s=1}^i
\frac{Kq^{-s+1}-K^{-1}q^{s-1}}
{q^s-q^{-s}}$, on a donc 
\begin{align}
\kappa'\kappa
&=
\frac{1}{2}\sum_{i=0}^{l-1}    (-1)^{i}     \left[\begin{array}{c}K \\i\end{array}\right] (q^{i} +   q^{ - i} K)
\kappa
=
\frac{1}{2}\sum_{i=0}^{l-1}    (-1)^{i}     \left[\begin{array}{c}K \\i\end{array}\right] (q^{i} +   q^{ - i})
\kappa
\\
\notag&=
\frac{1}{2}\sum_{i=0}^{l-1}    (-1)^{i} (q^{i} +   q^{ - i})
\prod_{s=1}^i
\frac{Kq^{-s+1}-K^{-1}q^{s-1}}
{q^s-q^{-s}} \kappa
\\
\notag&=
\frac{1}{2}\sum_{i=0}^{l-1}    (-1)^{i} (q^{i} +   q^{ - i})
\prod_{s=1}^i
\frac{q^{-s+1}-q^{s-1}}
{q^s-q^{-s}} \kappa
=\kappa,
\end{align}
comme voulu.

\subsection{}\label{}

Pour la suite, il est commode d'introduire, pour $n\in \mZ$,
 \begin{equation}\label{defkappaj}
 \kappa'_{n} =    \kappa_{2n}= \frac{1}{2l} \sum_{i=0}^{2l-1}q^{-ni}K^{i} \in U_{q}. 
 \end{equation}
On a clairement $\kappa'_n=\kappa'_m$ si     $n\equiv m\mod l$ et l'on \'etend donc la notation  $ \kappa'_{n}$ au cas o\`u $n \in \mZ/l$.
 
\begin{prop}\label{} 
On a, pour tout $n \in \mZ$, 
\begin{equation}\label{expresskappa}
\kappa_{2n} = \kappa'_n =      \frac{1}{2}\sum_{i=0}^{l-1}
(-1)^i\lbr K;-n
\\
i\rbr(q^i+q^{-i-n}K)\in u^0_\cB.
\end{equation}
\end{prop}

Comme $\kappa'_n$ ne d\'epend que de la r\'eduction modulo
$l$ de $n$, on peut   se limiter \`a \'etablir la formule lorsque $n$ est pair. Tout d'abord, on d\'eduit imm\'ediatement de la relation   $KF=q^{-2}FK$ \eqref{relgen}    la relation  de commutation
\begin{equation}\label{commutkappaF}
 F^{(n)}. \kappa'_{2n} =  \kappa.  F^{(n)}.
\end{equation}

On \'ecrit alors  
\begin{equation}
\begin{split}
F^{(n)}.\kappa'_{2n}
&=
\kappa .
F^{(n)} = (\frac{1}{2}. \{  \sum_{i=0}^{l-1}    (-1)^{i}     \left[\begin{array}{c}K \\i\end{array}\right] (q^{i} +   q^{ - i} K)  \}). F^{(n)}
\\
&=
\frac{1}{2}
\{
\sum_{i=0}^{l-1}
(-1)^i\lbr
K
\\
i\rbr
F^{(n)}
(q^i+q^{-i}q^{-2n}K)\}
\quad\text{en utilisant  (\ref{relgen}) }
\\
&=
F^{(n)}
\frac{1}{2}
\{
\sum_{i=0}^{l-1}
(-1)^i\lbr
K;-2n
\\
i\rbr
(q^i+q^{-i-2n}K)\}
\quad\text{par \cite[4.1.c]{L89}},
\end{split}
\end{equation}
on doit donc avoir 
$\kappa'_{2n}=\frac{1}{2}
\{
\sum_{i=0}^{l-1}
(-1)^i\lbr
K;-2n
\\
i\rbr
(q^i+q^{-i-2n}K)\}
$.


\subsection{}\label{} 
 On va aussi \'etablir l'int\'egralit\'e de $\kappa_{n}$ lorsque $n$ est impair. Pour cela, soit  \,\,$\tilde{}$ \,\,   l'automorphisme de  $U$ (\cite[4.6]{L89})     tel que  $E\mapsto-E, F\mapsto F$ and $K\mapsto-K$.    
 
\begin{lem} \label{kappabarren}
Pour tout $n\in {\Bbb{Z}}$, on a
\begin{equation}
\tilde{\kappa'_{n}} =\frac{1}{2}\sum_{i=0}^{l-1}
 \lbr K;-n
\\
i\rbr(q^i-q^{-i-n}K)\in u^0_\cB.
\end{equation}
\end{lem}
Il suffit d'appliquer la transformation de $K$ en $-K$ dans \eqref{expresskappa} et de remarquer que celle-ci transforme $ \lbr K;-n
\\
i\rbr$ en $ \lbr K;-n
\\
i\rbr$ pour $i$ pair et en  $- \lbr K;-n
\\
i\rbr$ pour $i$ impair.

\subsection{}\label{choixderacine}  
Dans la proposition suivante,  $\frac{n}{2}$ d\'esigne l'unique \'el\'ement de $\mZ/l$ dont le double est la classe de $n$ dans $\mZ/l$.

\begin{prop} 
\begin{itemize}
\item[{\rm (i)}] 

On a

\begin{equation}
 \kappa_n 
=
\begin{cases}
 \kappa'_{\frac{n}{2}}
&\text{si  $n$ est pair}
\\
\tilde\kappa'_{\frac{n}{2}}
&\text{si $n$ est impair}.
\end{cases}
\end{equation}
et par cons\'equent $\kappa_{n} \in u^0_\cB$ pour tout $n\in \mZ$ (ou tout $n\in \mZ/l$).

\item[{\rm (ii)}] 
On a
\begin{equation}
\widetilde{\kappa_n}
=
\kappa_{n+l}.
\end{equation}

 \end{itemize} 
\end{prop}

Pour (i), c'est une cons\'equence imm\'ediate   des d\'efinitions de $\kappa'_{n}$ \eqref{defkappaj} et de $\kappa_{n}$ \eqref{defkappa'j}.
Pour (ii),
si $n$ est pair,
\begin{equation}
\begin{split}
\widetilde{\kappa_n}
&=
\widetilde{\kappa'_{\frac{n}{2}}}
\quad\text{par  (i)}
\\
&=
\frac{1}{2}\sum_{i=0}^{
l-1}
(-1)^i\lbr -K;-\frac{n}{2}
\\
i\rbr(q^i-q^{-i-\frac{n}{2}}K)
\quad\text{gr\^ace \`a  \eqref{expresskappa}}
\\
&=
\frac{1}{2}\sum_{i=0}^{
l-1}
\lbr K;-\frac{n}{2}
\\
i\rbr(q^i-q^{-i-\frac{n}{2}}K)
=
\tilde\kappa'_{\frac{n}{2}}
=\kappa_{n+l}.
\end{split}
\end{equation}
Le cas $n$ impair se traite par un argument analogue.

\begin{cor}
\begin{itemize}


\item[{\rm (i)}] 
L'alg\`ebre $U_\cB$ admet une  d\'ecomposition
en somme directe 
\begin{equation}\label{decompU_cB}
U_\cB=\coprod_{n,m=0}^{2l-1}
\kappa_nU_\cB\kappa_m
\end{equation}
telle que 
$\Fr(\kappa_nU_\cB\kappa_m)\ne0$
si et seulement si 
$n=m=0$
et telle que 
$\phi: \bU_\cB\to U_\cB$ se  factorise \`a travers 
$\kappa_0U_\cB\kappa_0=
\kappa
U_\cB\kappa$.
\item[{\rm (ii)}] 
Posons $\kappa^+=\sum_{m\in[0,l[}\kappa_{2m}$
et
$\kappa^-=\sum_{m\in[0,l[}\kappa_{2m+1}$.
Les \'el\'ements 
$\kappa^+$ et $\kappa^-$ sont des idempotents centraux de  $U_\cB$
tels que
$\widetilde{\kappa^+}=\kappa^-$.
La  $\cB$-alg\`ebre $U_\cB$ admet une d\'ecomposition en somme d'id\'eaux
$U_\cB=(\kappa^+U_\cB\kappa^+)
\oplus
(\kappa^-U_\cB\kappa^-)=
(U_\cB\kappa^+)
\oplus
(U_\cB\kappa^-)=
(U_\cB\kappa^+)
\oplus
\widetilde{U_\cB\kappa^+}$
Tout $U_q$-module $M$ se d\'ecompose en  une somme directe
$M=(\kappa^+M)\oplus(\kappa^-M)$ telle  que
$K^l$ 
agisse sur  $\kappa^+
M$ (resp. $\kappa^-
M$)
par $1$ (resp. $-1$).
\end{itemize}
\end{cor}

Pour (ii), on remarquera simplement que pour tout
$r\in\bbN$ on a

\begin{align}
E^{(r)}\kappa_n
&=
E^{(r)}\frac{1}{2l}\sum_{i=0}^{2l-1}q^{\frac{-ni}{2}}K^i
=
\frac{1}{2l}\sum_{i=0}^{2l-1}q^{\frac{-ni}{2}-2ri}K^iE^{(r)}
\\
\notag&=
\frac{1}{2l}\sum_{i=0}^{2l-1}q^{\frac{-(n+4r)i}{2}}K^iE^{(r)}
=
\kappa_{n+4r}E^{(r)}.
\end{align}



De m\^eme, on prouve que
\begin{equation}
F^{(r)}\kappa_n=
\kappa_{n-4r}F^{(r)}.
\end{equation}

On a
\begin{equation} 
K\kappa_n=\frac{1}{2l}\sum_{i=0}^{2l-1}q^{\frac{-ni}{2}}K^{i+1}=q^{\frac{n}{2}}\kappa_n
=(-1)^nq^{\frac{(1-l)n}{2}}\kappa_n,
\end{equation}
et en  particulier,
$K^l\kappa_n=(-1)^n\kappa_n$.

\subsection{}\label{type}  

Soit 
$u_q
$
la $\bbQ(q)$-sous-alg\`ebre de
$U_q$ \eqref{defUq} engendr\'ee par
$E, F, K$. Tout $u_q$-module $M$ de dimension finie admet une d\'ecomposition en sous-espaces propres relativement \`a l'action de $K$ avec des valeurs propres des puissances de $q^{\frac{1}{2}}$.

\begin{defi}
On dit qu'un 
 $u_q$-module de dimension finie
$M$ est 
de type $1$ (resp. -1)
si et seulement si
$M=\kappa^{+}
M$ (resp. $M=\kappa^{-}
M$) .
\end{defi}

Cette d\'efinition est coh\'erente avec celle donn\'ee dans   \cite[4.6]{L89} et un $u_q$-module de dimension finie de type 1 n'est pas autre  chose qu'un $\kappa^+u_q\kappa^+$-module de dimension finie. Si  $M$ est un  
$u_q$-module 
de  type $-1$, on remarquera qu'en tordant l'action de
  $u_q$ par l'automorphisme 
 $\tilde\ $, 
on obtient un $u_q$-module  
de type  $1$
qu'on notera 
$\widetilde M$.

\subsection{}
Tout $u_q$-module 
ind\'ecomposable injectif/projectif est facteur direct d'un 
 $u_q\kappa_n$ avec $n\in[0,2l[$.
Rappelons (\cite{APW92}) maintenant quelques g\'en\'eralit\'es sur les 
 $u_q$-modules 
de   type 1.
Soit
$u_q^{\geq0}$
la sous-alg\`ebre de $u_q$
engendr\'ee par
$E$ et $K$. Pour $m \in \bbZ$, d\'esignons encore par $m$ le  
$u_q^{\geq0}$-module $\bbQ(q)$ tel que  $K\cdot1=q^m$ et $E\cdot1=0$,
et consid\'erons le $u_q$-module  $\bar\Delta_q(m)=u_q\otimes_{u_q^{\geq0}}m$. Le $u_q$-module  simple de plus haut poids $m$ (qui est aussi le module de t\^ete de  $\bar\Delta_q(m)$) sera not\'e $\bar L_q(m)$.   On a
$\bar L_q(m)=\bar L_q(m')$ si et seulement si $m\equiv m'\mod l$.
On note
$Q_q(m)$ la couverture projective (qui est aussi l'enveloppe injective)  de
$\bar L_q(m)$.
On a alors $Q_q(l-1)= \bar L_q(l-1)=
\bar\Delta_q(l-1)$ qui n'est autre que le module de  Steinberg   que nous noterons en abr\'eg\'e  
$\St_q$.
Chaque  $Q_q(m)$ 
(resp.
$\bar\Delta_q(m)$)
pour  $m\in[0,l-1[$ est, quant \`a lui, une extension non scind\'ee de 
$\bar\Delta_q(l-m-2)$ (resp.
$\bar L_q(l-m-2)$) par
$\bar\Delta_q(m)$
(resp.
$\bar L_q(m)$).


\begin{prop}\label{descrquant}
Soit
$n\in[0,2l[$.
\begin{itemize}
\item[{\rm (i)}] 
Si $n=2m$ est pair, on a
\begin{equation}
u_q\kappa_n
=
\coprod_{r\in\bbZ\cap([0,\frac{m}{2}[\cup[m,\frac{l+m-1}{2}])}
Q_q(2r-m)
\text{
avec
$Q_q(2r-m)=\St_q$ pour }
r=\begin{cases}
\frac{l+m-1}{2}
&\text{si $m$ pair}
\\
\frac{m-1}{2}
&\text{si $m$ impair}
\end{cases}.
\end{equation}
\item[{\rm (ii)}] Si $n$ est impair, on a
\begin{equation}
\widetilde{
u_q\kappa_n}
=
u_q\kappa_{n+l}.
\end{equation}
\end{itemize}
\end{prop}

Pour (i), nous aurons \`a utiliser l'involution $\Omega$ 
de  $\cB$-
alg\`{e}bre  d\'efinie par 
\begin{equation}\label{defomega}
\Omega(E)=F, \,\, \Omega(F)=E,\,\, \Omega(K)=K^{-1}.
\end{equation}
Examinons l'action de  $F^{(l-1)}$ sur les vecteurs primitifs
$E^{(l-1)}F^{(r)}\kappa_n$
de poids $2(l-1-r)+m$.
On a 
\begin{align}
F^{(l-1)}
&
E^{(l-1)}F^{(r)}\kappa_n
=
\sum_{i=0}^{l-1}
E^{(l-1-i)}
(-1)^i
\lbr
K;-i+2(l-1)-1
\\
i
\rbr
F^{(l-1-i)}F^{(r)}\kappa_n
\\
\notag&\hspace{6cm}
\text{en appliquant 
$\Omega$ \eqref{defomega}
\`a \cite[4.1.a]{L89}}
\\
\notag&=
\sum_{i=0}^{l-1}
E^{(l-1-i)}
(-1)^i
\lbr
K;-i+2l-3
\\
i
\rbr
\lbr
l-1-i+r
\\
r
\rbr
F^{(l-1-i+r)}\kappa_n
\\
\notag&=
\sum_{i=0}^{l-1}
E^{(l-1-i)}
(-1)^i
\lbr
l-1-i+r
\\
r
\rbr
F^{(l-1-i+r)}
\lbr
K;-i+2l-3-2(l-1-i+r)
\\
i
\rbr
\kappa_n
\\
\notag&\hspace{6cm}
\text{gr\^ace \`a  \cite[4.1.c]{L89}}
\\
\notag&=
\sum_{i=0}^{l-1}
E^{(l-1-i)}
F^{(l-1-i+r)}(-1)^i
\lbr
l-1-i+r
\\
r
\rbr
\lbr
K;i-2r-1
\\
i
\rbr
\kappa_n
\\
\notag&=
\sum_{i=0}^{l-1}
E^{(l-1-i)}
F^{(l-1-i+r)}(-1)^i
\lbr
l-1-i+r
\\
r
\rbr
\lbr
m+i-2r-1
\\
i
\rbr
\kappa_n
\\
\notag&=
\sum_{i=0}^{l-1}
E^{(l-1-i)}
F^{(l-1-i+r)}
\lbr
l-1-i+r
\\
r
\rbr
\lbr
2r-m
\\
i
\rbr
\kappa_n.
\end{align}
Si 
$0\leq
r<\frac{m}{2}$ ou si
$m\leq r\leq\frac{l+m-1}{2}$,
avec $i=r$
on a
$\lbr 
l-1-i+r
\\
r
\rbr
\lbr
2r-m
\\
i
\rbr
=
\lbr
l-1
\\
r
\rbr
\lbr
2r-m
\\
r
\rbr\ne0$,
et donc
$\coprod_{r\in\bbZ\cap([0,\frac{m}{2}[\cup[m,\frac{l+m-1}{2}])}
\bar\Delta_q(m-2-2r) \subseteq
u_q\kappa_n$.
Alors
$u_q\kappa_n=
\coprod_{r\in\bbZ\cap([0,\frac{m}{2}[\cup[m,\frac{l+m-1}{2}])}
Q_q(2r-m)$
pour des raisons de dimension
avec
$Q_q(2r-m)=\St_q$ pour
$r=\frac{l+m-1}{2}$
(resp. 
$r=\frac{m-1}{2}$)
si $m$ est pair
(resp. si $m$ est impair).

\section{Le cas  $\mathfrak{sl}_2$ : scindage du Frobenius}{}\label{ }

\begin{prop}\label{scindsl2} 
Posons, pour all\'eger,  $\kappa =\kappa_{0}$. L'application $
\cB$-lin\'{e}aire
\begin{equation}\label{phisl2}
\phi :   {\bU}_{\cB} \rightarrow 
U_{\cB}
\end{equation}
induite par multiplicativit\'e par,
\begin{equation}\label{deffi}
\phi(X^{(i)})= E^{(li)}.\kappa  \,\, ;  \,\,  \phi(Y^{(i)})= F^{(li)}.\kappa   \,\, ;  
\,\,  \phi(\small{\left(\begin{array}{c}H \\ i \end{array}\right)})=   \small{\left[\begin{array}{c}K \\li\end{array}\right]}. \kappa 
\end{equation}  
pour tout $i\geq 0$ est un homomorphisme (non-unif\`ere) d'alg\`{e}bres qui scinde l'application de Frobenius ${\rm{Fr}}$  \eqref{defFr}.  
\end{prop} 

Remarquons les \'egalit\'es \'evidentes  $\kappa. E^{(li)}= E^{(li)}.  \kappa$ ; $\kappa. F^{(li)}= F^{(li)}.  \kappa$ ;  $\small{\left[\begin{array}{c}K \\li\end{array}\right]}. \kappa  =  \kappa . \small{\left[\begin{array}{c}K \\li\end{array}\right]}$      pour tout  $i\geq0$ : dans la d\'{e}finition de $\phi$ \eqref{deffi}, on pourrait donc tout aussi bien  multiplier par $\kappa$ \`{a} gauche.

\subsection{}\label{}
 Si l'on veut formuler l'existence de cet homomorphisme en terme  d'homomorphisme  unif\`ere  d'alg\`ebres comme dans  Thm. \ref{teointrod} (ii), on peut proc\'eder comme suit. On remarque que l'on a une d\'ecomposition 
 $
U_\cB=\coprod_{i,j=0}^{2l-1}
\kappa_iU_\cB\kappa_j
$
de
$U_\cB$ en $\cB$-sous-modules
\eqref{decompU_cB}.
Comme ${\rm{Fr}}(\kappa_n)=\delta_{n0}$
pour tout $n\in[0,2l[$,
on obtient un diagramme commutatif d'homomorphismes de
$\cB$-alg\`ebres
\begin{equation}
\xymatrix{
U_\cB\ar[rr]^-{Fr}
\ar[d]
&&
 \bU_\cB
\\
U_\cB/\ker({\rm{Fr}})
\ar@{-->}[rru]_-\sim
\\
\kappa
U_\cB\kappa/\{
(\kappa
U_\cB\kappa)\cap\ker({\rm{Fr}})\}
\ar[u]^-\sim
&&
\kappa
U_\cB\kappa
\ar[ll]
\ar[uu]_-{{\rm{Fr}}\mid_{\kappa
U_\cB\kappa}}
}
\end{equation}
L'application 
$\phi$ va  alors \^etre un scindage (comme homomorphisme unif\`ere de $\cB$-alg\`ebres) de
${\rm{Fr}}|_{\kappa U_\cB\kappa}$.

\subsection{}\label{} 
Soit, comme dans l'introduction,  
$U^{\geq0}_{\cB}$
(resp.  
$U^{\leq0}_{\cB}$) la $\cB$-sous-alg\`{e}bre de   
$U_\cB$
engendr\'{e}e par  les $\{E^{(i)}, i\in \mN\}$ et les 
$\{K^{\pm 1}, \left[\begin{array}{c}K ;c \\m\end{array}\right], c\in \mZ, m\in \mN\}$ (resp. les $F^{(i)}$ et les $\{K^{\pm 1}, \left[\begin{array}{c}K ; c \\m\end{array}\right], c\in \mZ, m\in \mN$\}) et ${\bU}_{\cB}^{\geq0}$ (resp.  ${\bU}_{\cB}^{\leq0}$) la $
\cB$-sous-alg\`{e}bre de ${\bU}_{\cB}$ engendr\'{e}e par  les $\{X^{(i)}, i\in \mN\}$ et les $\{  \left(\begin{array}{c}H \\m\end{array}\right), m\in \mN\}$ 
(resp. les $\{Y^{(i)}, i\in \mN\}$ et les $\{  \left(\begin{array}{c}H \\m\end{array}\right), m\in \mN\}$).
Avec ces notations, il r\'{e}sulte alors de   \cite[Thm. 1.2]{kl}   que l'application (rappelons ici que $K^{l}$ est central dans $U_{\cB}$) 
\begin{equation}\label{descrFrob}
{\text{Fr}}' : {\bU}_{\cB}^{\geq0} \rightarrow  
U_\cB^{\geq0}/\{(K^{l}-1)(U_\cB^{\geq0}\otimes_\cB \bbQ(q))\cap
U_\cB^{\geq0}\}
\end{equation}
d\'{e}finie par 
\begin{equation}\label{descrFrob}
{\text{Fr}}'(X^{(i)})= E^{(l
i)}  \,\,,\,\,  {\text{Fr}}'(\left(\begin{array}{c}H \\i\end{array}\right) )=  \left[\begin{array}{c}K ;0 \\l
i\end{array}\right]  = \left[\begin{array}{c}K \\l
i\end{array}\right]         .
\end{equation}
et son analogue \'{e}vident   
$\bU_\cB^{\leq0}\rightarrow  
U_\cB^{\leq0}/\{(K^{l}-1)(U_\cB^{\leq0}\otimes_\cB \bbQ(q))\cap
U_\cB^{\leq0}\}$ sont des homomorphismes d'alg\`{e}bres.

 Comme la restriction de $\phi$ \`{a} $ {\bU}_{\cB}^{\geq0}$ et $ {\bU}_{\cB}^{\leq0} $ se factorise par ${\text{Fr}}'$, on voit donc, compte tenu de la remarque suivant la prop. \ref{scindsl2},  que cette restriction est un homomorphisme d'alg\`{e}bres. On utilise maintenant la relation (\cite[26.2]{Hum})
\begin{equation}\label{}
 X^{(b)} Y^{(a)} - Y^{(a)}  X^{(b)}   = \sum_{r=1}^{{\text{min}}(a,b)}Y^{(a-r)} \small{\left(\begin{array}{c}H-a-b+2r \\r\end{array}\right) }X^{(b-r)},
\end{equation}
pour tous $a\in \mN$ et $b\in \mN$. Il suffit donc  pour prouver   la 
prop. \ref{scindsl2}, de montrer qu'on a l'\'{e}galit\'{e} suivante dans  
$U_{\cB}$
\begin{equation}\label{compat1}
 \phi(X^{(b)})  \phi(Y^{(a)}) -  \phi(Y^{(a)})   \phi(X^{(b)})   = \sum_{r=1}^{{\text{min}}(a,b)}\phi(Y^{(a-r)}) \phi(\small{\left(\begin{array}{c}H-a-b+2r \\r\end{array}\right) })\phi(X^{(b-r)}).
 \end{equation}
 D'apr\`{e}s  (\cite{lus}, 6.5. (a2)), on a 
 \begin{equation}\label{}
 E^{(l b)}  F^{(l a)}  -  F^{(l a)}  E^{(l b)}   = \sum_{s=1}^{{\text{min}}(l a,l b)} F^{(l a-s)}  \small{ \left[\begin{array}{c}K ;  -l a-l b +2s\\s\end{array}\right] } E^{(l b-s)}.
  \end{equation} 
Utilisant la relation \eqref{commutkappaF},    on voit imm\'{e}diatement que  la v\'erification de \eqref{compat1} revient \`{a} \'{e}tablir que       pour tout $s$ tel que ${\text{min}}(la,lb) \geq s \geq 1$ et ne divisant pas $l$, on a 
   \begin{equation}\label{eqfondam}
  \kappa'_{2(l a+l b-s)}\small{ \left[\begin{array}{c}K ;  -l a-l b +2s\\s\end{array}\right] } =  \kappa'_{ -2s}\small{ \left[\begin{array}{c}K ;  2s-l a-l b\\s\end{array}\right] } = 0.
  \end{equation} 
On conclut en utilisant   le 
\begin{lem}\label{keylemma} 
Pour tous $n, m  \in\bbZ$ et tout $t\in\bbN$ non divisible par $l$, on a 
\begin{equation}\label{}
\kappa'_n\lbr
K;-n+lm
\\
t\rbr=0.
\end{equation}
\end{lem}

Supposons tout d'abord  $t<l$ et posons, pour $a \in \mZ$,  $\left[\begin{array}{c}a \\t\end{array}\right] =  \prod_{s=1}^t\frac{q^{a-s+1}-q^{-a+s-1}}{q^s-q^{-s}}\in \cB$. Cet \'el\'ement est nul si $a$ est un multiple de $l$. On a alors
\begin{equation}\label{casfondam}
\begin{split}
\kappa'_n\lbr
K;-n+lm
\\
t\rbr
&=
\kappa_n\prod_{s=1}^t
\frac{Kq^{-n+lm-s+1}-K^{-1}q^{n-lm+s-1}}
{q^s-q^{-s}}
\\
&=
\prod_{s=1}^t
\frac{q^{lm-s+1}-q^{-lm+s-1}}
{q^s-q^{-s}}
\quad\text{gr\^ace \`a \ref{weightofkappa}}
\\
&=
\lbr
lm
\\
t\rbr
\\
&=
0.
\end{split}
\end{equation}
Plus g\'en\'eralement, \'ecrivons   
$t=t_0+t_1l$
avec $t_0\in[0,l[$ et $t_1\in\bbZ$.
On remarque alors qu'on dispose d'un analogue de  
\cite[g3]{lus2}, \`a savoir que pour tous $k,k'\in\bbN$, on a 
\begin{equation}
\lbr
K;m
\\
k
\rbr
\lbr
K;m-k
\\
k'
\rbr
=
\lbr
k+k'
\\
k
\rbr
\lbr
K;m
\\
k+k'
\rbr.
\end{equation}
Faute de r\'ef\'erence pour cette \'egalit\'e, signalons qu'il suffit de la v\'erifier dans    $U_{\bbQ(v)}$ :
\begin{equation}
\begin{split}
\lbr
K;m
\\
k
\rbr
\lbr
K;m-k
\\
k'
\rbr
&=
\prod_{s=1}^k\frac{Kv^{m-s+1}-K^{-1}v^{-m+s-1}}{v^s-v^{-s}}\prod_{s=1}^{k'}\frac{Kv^{m-k-s+1}-K^{-1}v^{-m+k+s-1}}{v^s-v^{-s}}
\\
&=
\prod_{s=1}^{k+k'}\frac{Kv^{m-s+1}-K^{-1}v^{-m+s-1}}{v^s-v^{-s}}\frac{\prod_{s=1}^{k+k'}(v^{s}-v^{-s})}{\prod_{s=1}^{k}(v^s-v^{-s})
\prod_{s=1}^{k'}(v^s-v^{-s})}
\\
&=
\lbr
K;m
\\
k+k'
\rbr
\lbr
k+k'
\\
k
\rbr.
\end{split}
\end{equation}
On en d\'eduit alors
\begin{equation}
\begin{split}
\kappa'_n\lbr
K;-n+lm
\\
t
\rbr
&=
\kappa'_n
\lbr
t
\\
t_1l
\rbr
\lbr
K;-n+lm
\\
t
\rbr
=
\lbr
K;-n+lm
\\
t_1l
\rbr
\kappa'_n
\lbr
K;-n+lm-lt_1
\\
t_0
\rbr
\\
&=
0
\quad\text{gr\^ace \`a (\ref{casfondam})}.
\end{split}
\end{equation}

\begin{rema}\label{}
Si l'on utilise le lemme \ref{keylemma}, on peut se dispenser d'invoquer  \cite[Thm. 1.2]{kl} pour voir que $\mathrm{Fr}'$ s'\'etend en une application multiplicative  $\phi:
\bU_\cB^{\geq0}\to
U_\cB^{\leq0}$ and $\phi:\bU_\cB^{\leq0}\to
U_\cB^{\leq0}$.
\end{rema}

\section{Retour sur la th\'eorie modulaire}\label{ }

\subsection{}\label{distrib}
Soit $p$ un nombre premier impair. La th\'eorie que nous qualifions ici   en abr\'eg\'e de \emph{modulaire}  r\'ef\`ere    toujours  au cas o\`u $l=p$ et  o\`u, \'eventuellement,  l'on ``r\'eduit'' modulo $p$ la th\'eorie quantique en un sens qu'on va pr\'eciser tout de suite \eqref{famillekappan}.  Soient alors  $G$ le $\bbF_{p}$-groupe alg\'ebrique $SL_{2}$ des matrices carr\'es d'ordre 2 et de d\'eterminant 1 et $T$ le tore maximal form\'e des matrices diagonales. Pour $r\in\bbN^+$, soient  $G_r$ (resp. $T_r$) le $r$-i\`eme noyau de Frobenius de  $G$ (resp. $T$)   et $\Dist(G)$ (resp. $\Dist(G_r)$, $\Dist(T)$, $\Dist(T_r)$) l'alg\`ebre des distributions de $G$ (resp. $G_r$, $T$, $T_r$).  On a, dans $\Dist(G)$, les relations suivantes pour tous $a,a',b,c,c' \in \mN$,
\begin{equation}\label{reldefDist1}
 X^{(a)}X^{(a')}=
\binom{a+a'}{a}X^{(a+a')},
\quad
Y^{(c)}Y^{(c')}=
\binom{c+c'}{c}Y^{(c+c')},
\end{equation}
\begin{equation}\label{reldefDist2}
X^{(a)}\binom{H}{b}=
\binom{H-2a}{b}X^{(a)},
\quad
Y^{(c)}\binom{H}{b}=
\binom{H+2c}{b}Y^{(c)}
\quad\text{\cite[Lem. 26.3.D]{Hum}},
\end{equation}
\begin{equation}\label{reldefDist3}
X^{(a)}Y^{(c)}=
\sum_{i=0}^{\min\{a,c\}}
Y^{(c-i)}\binom{H+2i-a-c}{i}X^{(a-i)}
\quad\text{\cite[Lem. 26.2]{Hum}}.
 \end{equation}

\subsection{}\label{Fr'}
Conservons les notations de \ref{distrib}. On dispose sur $\Dist(T)$ de l'endomorphisme (dit de Frobenius) not\'e $\Dist(\Fr)$ (ou parfois simplement $\Fr$) de ${\Bbb{F}}_{p}$-alg\`ebres tel que pour tout $m\in\bbN$
\begin{equation}
\binom{H}{m}\mapsto
\begin{cases}
\binom{H}{\frac{m}{p}}
&\text{si $p|m$}
\\
0
&\text{sinon}.
\end{cases}
\end{equation}
et de son scindage \'evident $\Fr'$ d\'efini par   
\begin{equation}\label{defFrprime}
\Fr'(\binom{H}{m}) =
\binom{H}{pm}
\end{equation} 
qui est un endomorphisme de ${\Bbb{F}}_{p}$-alg\`ebres comme on le v\'erifie ais\'ement. Il r\'esulte des propri\'et\'es de l'application   $\Fr'^{\pm}$ mentionn\'ee en \ref{intro3} que cette notation est coh\'erente  avec \eqref{descrFrob}.
 
\subsection{}\label{mun}
Nous aurons besoin plus bas (\ref{famillekappan}), pour    $n\in\bbZ$, de  
\begin{equation}\label{defmun}
\mu_n=\sum_{i=0}^{p-1}(-1)^i\binom{H-n}{i}\in\Dist(T_1).
\end{equation} 
 Ces \'el\'ements sont ceux qui \'etaient not\'es $\Delta_{T,n}$ dans \cite[3.1]{g}. Mentionnons \`a cette occasion que   \cite[Prop. 3.1.6]{g}  contient (au moins) une erreur typographique : avec les notations adopt\'ees ici, c'est   
$X^{(n)}\mu_m=
\mu_{m+2n}X^{(n)}$
et
$Y^{(n)}\mu_m=
\mu_{m-2n}Y^{(n)}$ qu'il faut lire.
Rappelons aussi  
(\cite[3.1.5]{g}) que
$\mu_n=\mu_m$
si et seulement si
$n\equiv m\mod p$. 
\subsection{}
 La relation avec  les  
$\kappa'_n  \in  U_q$,  $n\in\bbZ$
consid\'er\'es   en  \eqref{defkappaj} est donn\'ee par la

\begin{prop}\label{famillekappan}
Si $l=p$ est premier,
via l'homomorphisme  canonique
\begin{equation}\label{quantversmod}
U_\cB\to\Dist(G)
\end{equation}
 induit  par la projection canonique
$\cB\to\bbF_p$, ($q\mapsto1$) et tel  que $K\mapsto1$, l'image de $\kappa'_n$ est \'egale \`a $\mu_{n}$.
\end{prop}
Cela d\'ecoule aussi de \eqref{commutkappaF} jointe \`a 
\cite[Prop. 3.1.6]{g}.

 \subsection{}\label{relducasmod}
 On avait introduit dans la th\'eorie modulaire \cite[Thm. 1.3]{gk} un scindage (non unif\`ere) not\'e aussi (le contexte enlevant   tout risque de confusion avec \eqref{defscind} dans les notations)
 \begin{equation}\label{scindmod}
 \phi : \Dist(G) \to \Dist(G)
 \end{equation}
 de l'application canonique de Frobenius 
  \begin{equation}\label{Frobmod}
  \Dist(\Fr) = \Fr : \Dist(G) \to \Dist(G).
 \end{equation}\label{comparascind}
 Rappelons que cet endomorphisme $\phi$, restreint \`a $\Dist(T)$ s'obtient simplement en composant $\Fr'$ \eqref{defFrprime} et la multiplication par $\mu_{0}$ \eqref{defmun}.
 Il r\'esulte imm\'ediatement de \eqref{famillekappan}   sp\'ecialis\'e au cas $n=0$ qu'on a 
\begin{cor}\label{relev}
Pour $l=p$, l'application $\phi$ quantique \eqref{defscind} rel\`eve, via \eqref{quantversmod},  l'application $\phi$  modulaire \eqref{scindmod}.
\end{cor}
\subsection{} 
On revient \`a la situation de \ref{distrib}-\ref{mun}.
\begin{prop}\label{decomp1modulaire} 
On a 
\begin{equation}
\sum_{n=0}^{p-1}\mu_n=1\in\Dist(T_1)
\end{equation}
 et cette d\'ecomposition est une d\'ecomposition en idempotents deux \`a deux orthogonaux.
 Les  $\mu_n$, $n\in[0,p[$,
forment une base orthogonale de   
$\Dist(T_1)$. 
Par rapport \`a la base standard de $\Dist(T_1)$ form\'ee des $\binom{H}{i}$, $i\in[0,p[$, on a les formules de changement de base pour tout $n\in [0,p[$
\begin{equation}  
\mu_n=
\sum_{i=0}^{p-1}\binom{p-1-n}{p-1-i}\binom{H}{i} ,
\end{equation}
\begin{equation}  
\binom{H}{n}
=
\sum_{i=n}^{p-1}
\binom{i}{n}\mu_i.
\end{equation}
\end{prop}

Le fait que  chaque $\mu_n$ soit un idempotent r\'esulte de \cite[3.1.4]{g} et de \cite[3.1.3]{g} qu'on ait  
\begin{equation} 
\begin{split}
\mu_n
&=
\binom{H-n-1}{p-1}=
\frac{(H-n-1)(H-n-2)\dots(H-n-p+1)}{(p-1)!}
\\
&=
-(H-n-1)(H-n-2)\dots(H-n-p+1).
\end{split}
\end{equation}
Soient $n<m$,    $(n, m) \in [0,p[^{2}$.
Pour v\'erifier que $\mu_m\mu_n=0$, il suffit de voir dans  $\Dist(T_1)$ que
\begin{multline}\label{pdtH}
(H-m-1)(H-m-2)\dots(H-m-p+1)
\\
(H-n-1)(H-n-2)\dots(H-n-p+1)
=0.
\end{multline}
Comme
$-m-p+1\leq-n$
, on remarque alors que le membre de gauche de 
\eqref{pdtH} contient le produit
$H(H-1)\dots(H-p+1)=\binom{H}{p}p!=0$.
 
D'autre part, d'apr\`es   \cite[Cor. 3.1.5, Cor. 3.1.3 et Prop. 3.1.1]{g} on a 
$\mu_n
=\mu_{n-p}
=\binom{H+p-n-1}{p-1}=
\sum_{i=0}^{p-1}\binom{p-n-1}{p-i-1}\binom{H}{i}$.
Comme la matrice 
$[\!(\binom{p-n-1}{p-i-1})\!]_{(n,i)\in [0,p[^{2}}$ 
 est unipotente,
les $\mu_n$, $n\in[0,p[$, forment une base de $\Dist(T_1)$.
Ecrivant 
$\binom{H}{n}=\sum_{i=0}^{p-1}
c_i\mu_i=
\sum_{i=0}^{p-1}
c_i\binom{H-i-1}{p-1}$ et appliquant aux deux c\^ot\'es 
les caract\`eres $\bar\chi_{j}$,
$j\in\bbZ$,
d\'efinis pour tout $i\in\bbN$ par
\begin{equation}\label{chi}
\bar\chi_{j}(\binom{H}{i})=\binom{j}{i},
\end{equation}  
on obtient
$[\!(\binom{p-n-1}{p-i-1})\!]^{-1}=
[\!(\binom{i}{n})\!]$.
En particulier,
$1=\sum_{i=0}^{p-1}\mu_i$.

\subsection{}\label{}  
Nous \'etablirons  plus loin \eqref{descrmodul} l'analogue modulaire de \ref{descrquant} mais l'on va tout d'abord  introduire la g\'en\'eralisation suivante de  \eqref{decomp1modulaire} qui nous sera utile plus bas. Posons, pour tous $r\in\bbN^+$ et $n\in\bbZ$,
\begin{equation}\label{defmunr}
\mu_n^{(r)}=\sum_{i=0}^{p^r-1}(-1)^i\binom{H-n}{i} \in \Dist(T_r),
\end{equation} 
si bien que 
$\mu_n^{(1)}$ d\'esigne le  $\mu_n$ introduit en \eqref{defmun}.

\begin{prop}\label{munsup}
Soit $r\in\bbN^+$. 
\begin{itemize}
\item[{\rm (i)}]
Pour tout $n\in\bbZ$, on a
$\mu_n^{(r)}=\binom{H-n-1}{p^r-1}=\sum_{i=0}^{p^r-1}\binom{-1}{p^r-1-i}\binom{H-n}{i} \in \Dist(T_r)$.
\item[{\rm (ii)}]
Pour tous $n, m\in\bbZ$, on a
$\mu_n^{(r)}=\mu_m^{(r)}$
si et seulement si $n\equiv m\mod p^r$.
\item[{\rm (iii)}] Soient $\phi$   le scindage de Frobenius modulaire \eqref{scindmod} et $\Fr'$ l'homomorphisme \eqref{defFrprime}. Pour tous 
$m\in[0,p[$ 
et $n\in\bbZ$,
$\mu_{m+np}^{(r+1)}=\mu_m{\Fr'}(\mu_n^{(r)})$,
et donc
$\phi(\mu_n^{(r)})=\mu_{np}^{(r+1)}$.
\item[{\rm (iv)}]
Les $\mu_n^{(r)}$ pour $n\in[0,p^r[$,
fournissent une d\'ecomposition de   1 en idempotents orthogonaux deux \`a deux. En particulier,   $\mu_0^{(r)}$ est l'unique mesure involutive invariante de 
$\Dist(T_r)$,
i.e.,
pour tout $\mu\in\Dist(T_r)$,
$\mu\mu_0^{(r)}=\varepsilon(\mu)\mu_0^{(r)}$
avec
$\varepsilon=\bar\chi_0$ \eqref{chi}  la co\"unit\'e de $\Dist(T_r)$.



\item[{\rm (v)}]Les $\mu_n^{(r)}$,  pour $n\in[0,p^r[$,
forment   une base    orthogonale de
$\Dist(T_r)$. Par rapport \`a la base standard de $\Dist(T_r)$ form\'ee des $\binom{H}{i}$, $i\in[0,p^{r}[$, on a les formules de changement de base pour tout $n\in [0,p^{r}[$ : 
 \begin{equation}\label{expmunr}
\mu_n^{(r)}=\mu_{n-p^r}^{(r)}= 
\sum_{i=0}^{p^r-1}\binom{p^r-1-n}{p^r-1-i}\binom{H}{i},
\end{equation}
 \begin{equation}
\binom{H}{n}
=
\sum_{i=0}^{p^r-1}
\binom{i}{n}\mu_i^{(r)}
=
\sum_{i=n}^{p^r-1}
\binom{i}{n}\mu_i^{(r)}.
\end{equation}
Plus g\'en\'eralement, pour tous $m\in\bbZ$ et $n\in[0,p^r[$,
\begin{equation}
\binom{H-m}{n}
=
\sum_{i=0}^{p^r-1}
\binom{i-m}{n}\mu_i^{(r)}
\end{equation}
\item[{\rm (vi)}]
Pour tout
$n\in[0,p^r[$, on a 
$\Dist(\Fr^r)(\mu_n^{(r)})=
\Dist(\Fr)^r(\mu_n^{(r)})=\delta_{n0}$.
\end{itemize}
\end{prop}
 
Pour (i), la seconde \'egalit\'e r\'esulte de  
 \cite[Cor. 3.1.2]{g}.
Comme on a aussi
 \begin{equation}
\binom{-1}{p^r-1-i}=
\frac{(-1)(-2)\dots(-1-(p^r-1-i)+1)}{(p^r-1-i)!}
=
(-1)^{p^r-1-i}\frac{(p^r-1-i)!}{(p^r-1-i)!}
=
(-1)^i,
\end{equation}
l'assertion  en d\'ecoule.
Pour (ii), on prouve d'abord le sens ``si" pour lequel il suffit de montrer que
 $\mu_{n-p^r}^{(r)}=\mu_n^{(r)}$. Or, 
\begin{equation}
\begin{split}
\mu_{n-p^r}^{(r)}
&=
\binom{H+p^r-n-1}{p^r-1}
\quad\text{par (i)}
\\
&=
\sum_{i=0}^{p^r-1}\binom{p^r}{p^r-1-i}\binom{H-n-1}{i}
\quad\text{gr\^ace \`a  \cite[Cor. 3.1.2]{g} de nouveau}
\\
&=
\binom{H-n-1}{p^r-1}
\\
&=
\mu_n^{(r)}
\quad\text{par (i)}.
\end{split}
\end{equation}

Pour (iii), 
soient $m\in[0,p[$ et $n\in\bbZ$.
Gr\^ace \`a ce qu'on vient juste de d\'emontrer, on peut supposer   
$n\in[0,p^r[$.
On a 
\begin{equation}
\begin{split}
\mu_{m+np}^{(r+1)}
&=
\mu_{m+np-p^{r+1}}^{(r+1)}
\quad\text{partie  ``si" de (ii)
}
\\
&=
\binom{H+p^{r+1}-m-np-1}{p^{r+1}-1}
\quad\text{par (i)}
\\
&=
\sum_{i=0}^{p^{r+1}-1}
\binom{p^{r+1}-m-np-1}{p^{r+1}-1-i}
\binom{H}{i}
\quad\text{par \cite[Prop. 3.1.1]{g}}
\\
&=
\sum_{i=0}^{p^{r+1}-1}
\binom{p(p^r-n-1)+p-m-1}{p(p^r-i_1-1)+p-i_0-1}
\binom{H}{i_0+i_1p}
\quad\text{si l'on \'ecrit
$i=i_0+i_1p$}
\\
&\hspace{2cm}
\text{avec $i_0\in[0,p[$ et  $i_1\in\bbN$}
\\
&=
\sum_{i=0}^{p^{r+1}-1}
\binom{p^r-n-1}{p^r-i_1-1}\binom{p-m-1}{p-i_0-1}
\binom{H}{i_0}\binom{H}{i_1p}
\\
&=
\mu_m\Fr'(\mu_n^{(r)}).
\end{split}
\end{equation}

La premi\`ere assertion de 
(iv) d\'ecoule maintenant de   \eqref{decomp1modulaire} 
gr\^ace \`a  (iii) car $\Fr'$ est un homomorphisme de  $\bbF_p$-alg\`ebres.
La partie ``seulement si " de (ii) 
en d\'ecoule \'egalement. Quant aux formules de changement de base (v), elles se prouvent 
comme dans la proposition \eqref{decomp1modulaire} \`a l'aide des $\bar\chi_{ j}$.
L'assertion (vi) d\'ecoule imm\'ediatement de (v).

\subsection{}
A l'aide  de cette base orthogonale   de $\Dist(T_r)$ form\'ee par les 
$\mu_n^{(r)}$ \eqref{defmunr}, on peut  voir que les $X^{(a)}\mu_{b}^{(r)}Y^{(c)}$ forment, pour $a, b, c\in[0,p^r[$,  une base de $\Dist(G_r)$. Outre  \eqref{reldefDist1}, les relations  entre ces \'el\'ements (correspondantes \`a \eqref{reldefDist2}-\eqref{reldefDist3}) d\'efinissant $\Dist(G_r)$  s'expriment comme suit.    

\begin{prop}\label{commutamuXY}
Soient $a, b, c\in[0,p^r[$.
\begin{itemize}
\item[{\rm (i)}]

$X^{(a)}\mu_b^{(r)}=
\mu_{b+2a}^{(r)}X^{(a)}$,
\quad
$Y^{(c)}\mu_b^{(r)}=
\mu_{b-2c}^{(r)}Y^{(c)}$.

\item[{\rm (ii)}]
$X^{(a)}Y^{(c)}
=
\sum_{i=0}^{\min\{a,c\}}
\sum_{j=0}^{p^r-1}
\binom{j-a-c+2i}{i}
Y^{(c-i)}\mu_j^{(r)}X^{(a-i)}$.

\end{itemize}
\end{prop}
 
(i) Combinant  \eqref{expmunr} et \eqref{reldefDist2}, 
on  a en effet pour tout $n\in\bbZ$
\begin{equation}
\begin{split}
X^{(a)}\mu_n^{(r)}
&=
X^{(a)}\binom{H-n-1}{p^r-1}
 \quad\text{par \ref{munsup} (i)}
\\
&=
\binom{H-2a-n-1}{p^r-1}X^{(a)}
 \quad\text{par \eqref{reldefDist2}}
\\
&=
\mu_{n+2a}^{(r)}
X^{(a)}
 \quad\text{par \ref{munsup} (i)}.
\end{split}
\end{equation}
Par un calcul analogue (ou en utilisant l'involution de Chevalley),
\begin{equation}
Y^{(c)}\mu_n^{(r)}=
\mu_{n-2c}^{(r)}Y^{(c)}.
\end{equation}

(ii)  d\'ecoule de  \eqref{reldefDist3} et de la formule de changement de base \eqref{expmunr}.

\begin{cor}\label{decompDistGr}
Soit $r\in\bbN^+$. 
\begin{itemize}
\item[{\rm (i)}]
Si  $1 \leq s\leq r$,
$\Dist(G_r)$ admet une d\'ecomposition en somme directe 
\begin{equation}
\Dist(G_r)=\coprod_{n,m=0}^{p^s-1}\mu_n^{(s)}\Dist(G_r)\mu_m^{(s)}
\end{equation}
et chaque  $\mu_n^{(s)}\Dist(G_r)\mu_m^{(s)}$ est un  $\bbF_p$-espace vectoriel de base
$X^{(a)}Y^{(c)}\mu^{(r)}_{kp^s+m}$ avec $k\in[0,p^{r-s}[$,
$a, c\in[0,p^r[$ tels que
$n+2c\equiv m+2a\mod p^s$. 
\item[{\rm (ii)}]
$\Dist(G)$ admet une d\'ecomposition en somme directe 
\begin{equation}
\Dist(G)=\coprod_{n,m=0}^{p^r-1}\mu_n^{(r)}\Dist(G)\mu_m^{(r)}
\end{equation}
et chaque  $\mu_n^{(r)}\Dist(G_{r+s})\mu_m^{(r)}$,
$s\in\bbN^+$,
est un  $\bbF_p$-espace vectoriel de base
$X^{(a)}Y^{(c)}\mu^{(r+s)}_{kp^r+m}$,
  avec
$k\in[0,p^s[$, $a,c\in[0,p^{r+s}[$
tels que
$n+2c\equiv m+2a\mod p^r$. 
\end{itemize}
\end{cor}
 
La proposition \ref{munsup} (iv) permet d'\'ecrire
\begin{equation}
\begin{split}
\Dist(T_r)\mu_m^{(s)}
&=\coprod_{k=0}^{p^r-1}\bbF_p\mu_k^{(r)}\mu_m^{(s)}
\\
&=
\coprod_{k=0}^{p^r-1}\bbF_p\mu_{k_0}^{(s)}(\Fr')^s(\mu_{k_1}^{(r-s)})\mu_m^{(s)}
\quad\text{avec
$k=k_0+k_1p^s$
et
$k_0\in[0,p^s[$,
 $k_1\in\bbN$}
\\
&=
\coprod_{k=0}^{p^{r-s}-1}\bbF_p\mu_{m}^{(s)}(\Fr')^s(\mu_{k}^{(r-s)})
=
\coprod_{k=0}^{p^{r-s}-1}\bbF_p\mu_{m+kp^{s}}^{(r)}.
\end{split}
\end{equation}
Et donc
\begin{equation}
\begin{split}
\mu_n^{(s)}X^{(a)}Y^{(c)}\mu_{m+kp^{s}}^{(r)}
&=
X^{(a)}Y^{(c)}\mu_{n-2a+2c}^{(s)}\mu_{m+kp^{s}}^{(r)}
\quad\text{par \ref{commutamuXY}}
\\
&=
\begin{cases}
X^{(a)}Y^{(c)}\mu^{(r)}_{m+kp^{s}}
&\text{si $n-2a+2c\equiv m\mod p^s$}
\\
0
&\text{sinon}.
\end{cases}
\end{split}
\end{equation}

\subsection{}\label{}  

Les arguments utilis\'es pour prouver le corollaire  \ref{decompDistGr}, sp\'ecialis\'es au cas $r=1$,  donne l'existence d'une 
d\'ecomposition
$\Dist(G_1)=\coprod_{n=0}^{p-1}\Dist(G_1)\mu_n$
de $\Dist(G_1)$ en somme de $\Dist(G_1)$-modules \`a la fois projectifs et injectifs que nous allons maintenant identifier.
Notons  $\bar L(m)$, $m\in\bbZ$ les  $\Dist(G_1)$-modules simples de plus haut poids  
  $m$ avec
$\bar L(m)=\bar L(m')$ si et seulement si $m\equiv m'\mod p$.
Soit 
$Q(m)$ la couverture projective (qui est aussi l'enveloppe injective)  de
$\bar L(m)$.
Soit
$\Dist(B_1^+)$
la sous-alg\`ebre de $\Dist(G_1)$
engendr\'ee par
$X$ et $H$. Pour $m \in \bbF_p$, d\'esignons encore par $m$ le  
$\Dist(B^+_1)$-module $\bbF_p$ tel que  $H\cdot1=m$ et $X\cdot1=0$,
et posons  alors $\bar\Delta(m)=\Dist(G_1)\otimes_{\Dist(B_1^+)}m$.
On a $Q(p-1)=\bar L(p-1)=\bar 
\Delta(p-1)$ qui n'est autre que le module de  Steinberg   que nous noterons en abr\'eg\'e  
$\St$.
Chaque  $Q(m)$ pour  $m\in[0,p-1[$ est, quant \`a lui, une extension non scind\'ee de 
$\bar\Delta(p-m-2)$ par
$\bar\Delta(m)$.
Le m\^eme argument que dans la proposition \ref{descrquant} fournit la proposition suivante.

\begin{prop}\label{descrmodul}
\begin{itemize}
\item[{\rm (i)}] 
Si $p=2$,
$\Dist(G_1)\mu_0=Q(0)$
alors que
$\Dist(G_1)\mu_1=\St\oplus\St$.

\item[{\rm (ii)}] Si $p$ est impair, on a,  pour tout $n\in[0,p[$ pair
\begin{align*}
\Dist(G_1)\mu_n
&=
\coprod_{m\in\bbZ\cap([0,\frac{n-2}{2}]\cup
[n,\frac{p+n-1}{2}])}Q(2m-n)
\quad\text{avec
$Q(2m-n)=\St$ pour $m=\frac{p+n-1}{2}$}.
\end{align*}
\item[{\rm (iii)}] Si $p$ est impair, on a, pour tout  $n\in[0,p[$ impair,
\begin{align*}
\Dist(G_1)\mu_n
&=
\coprod_{m\in\bbZ\cap([0,\frac{n-1}{2}]\cup
[n,\frac{p+n-2}{2}])}Q(2m-n)
\quad\text{avec
$Q(2m-n)=\St$ pour $m=\frac{n-1}{2}$}.
\end{align*}

\end{itemize}
\end{prop}

\subsection{}
Il r\'esulte de \ref{descrmodul} qu'aucun des  
$\Dist(G_1)\mu_n$,
$n\in\bbZ$  n'est un  prog\'en\'erateur pour $\Dist(G_1)$.
 Il en va  de m\^eme   avec
$\Dist(G_2)\mu_0$ pour $\Dist(G_2)$ comme l'on peut d\'ej\`a le v\'erifier   pour $p=2$. En effet,
dans le cas contraire, $\St_2=\bar\Delta_2(p^2-1)=
\bar\Delta_2(3)=\Dist(G_2)\otimes_{\Dist(B^+_2)}3$
devrait appara\^itre dans  une d\'ecomposition de $\Dist(G_2)\mu_0$. Raisonnons alors dans la   base de $\Dist(G_2)\mu_0$
  fournie par les $\{
X^{(a)}Y^{(b)}\mu_{cp}^{(2)}\mid
a, b\in[0,p^2[, c\in[0,p[\}$ avec $X^{(a)}Y^{(b)}\mu_{cp}^{(2)}$ de poids 
$2(a-b)+pc$. Si ${\text{Hom}}_{\Dist(G_2)}(\St_2, \Dist(G_2)\mu_0)\ne0$, l'image de $1\otimes1$
dans
$\Dist(G_2)\mu_0$ doit \^etre annul\'ee par  
 $\Dist^+(U_2^+)=\coprod_{n\in]0,p^2[}{\Bbb{F}}_{p}
X^{(n)}$. Cette image est une combinaison lin\'eaire de $X^{(a)}Y^{(b)}\mu_{cp}^{(2)}$ comme ci-dessus avec 
$a=p^2-1$ ; or le poids de $X^{(p^2-1)}Y^{(b)}\mu_{cp}^{(2)}$ est   $2(p^2-1-b)+pc$ qui n'est certainement pas \'egal \`a 3, le poids de $1\otimes1$ ; d'o\`u la contradiction.

\subsection{}

Pour tout nombre premier  $l=p$ impair, notons   $\widehat\cB$ le compl\'et\'e de  $\cB$ relativement \`a l'id\'eal  $(q-1)$.
Tout $\Dist(G_1)$-module ind\'ecomposable injectif/projectif est   facteur direct d'un 
$\Dist(G_1)\mu_n$ avec $n\in[0,p[$.
Si maintenant $u_\cB$ d\'esigne  la sous-alg\`ebre de  
$U_\cB$ engendr\'ee par
$E, F, K$ et si l'on pose $u_{\widehat\cB}=u_\cB\otimes_\cB\widehat\cB$,
alors $u_{\widehat\cB}\kappa'_{2n}$ est un rel\`evement de
$\Dist(G_1)\mu_n$.
Il r\'esulte de 
\cite[5.4]{APW92} que les $\Dist(G_1)
$-modules ind\'ecomposables 
injectifs/projectifs se rel\`event en un facteur   de 
$u_{\widehat\cB}\kappa'_{2n}$ 
pour $n$ convenable
(cf. \cite[5.6]{APW92}).
En sens inverse,   tout $u_q\otimes_\cB{\text{Frac}}(\widehat\cB)$-module int\'egrable ind\'ecomposable  
projectif est   facteur direct d'un $u_q\kappa'_n\otimes_\cB{\text{Frac}}(\widehat\cB)$. 
Ces arguments s'\'etendent au cadre plus g\'en\'eral dans lequel on   d\'emontrera le th\'eor\`eme \ref{teointrod} plus bas.

 \subsection{}

On peut aussi d\'eduire de  la simple existence des idempotents \eqref{defmunr} et de leurs propri\'et\'es une nouvelle preuve de l'existence du scindage de 
 $\Dist(\Fr)$ sur $\Dist(G)$ 
de  \cite{g}.  Commen\c{c}ons par le
 
\begin{lem} \label{lemmu0Dist}
L'ensemble 
$\sum_{a,b,c\in\bbN}
\sum_{n\in[0,p^{b}[}
\Bbb{F}_{p}
X^{(pa)}Y^{(pc)}\mu_{np^{b}}^{(1+b)}$
est une sous-$\Bbb{F}_{p}$-alg\`ebre de 
$\mu_0\Dist(G)\mu_0$.
\end{lem}

Par \ref{decompDistGr},
chaque
$X^{(pa)}Y^{(pc)}\mu_{np^{b}}^{(1+b)}$,
$a,b,c\in\bbN$,
$n\in[0,p^b[$, appartient \`a
$\mu_0\Dist(G)\mu_0$.
Pour tous $r,s,t\in\bbN$,
$m\in[0,p^s[$, prenant
$d>\max\{ b,s\}$
tel que
$p^d>\max\{pc,pr\}$,
on a, par \ref{commutamuXY}
\begin{align}
X^{(pa)}
&
Y^{(pc)}\mu_{np^{b}}^{(1+b)}
X^{(pr)}Y^{(pt)}\mu_{mp^{s}}^{(1+s)}
=
X^{(pa)}Y^{(pc)}X^{(pr)}Y^{(pt)}\mu_{np^{b}-2pr+2pt}^{(1+b)}
\mu_{mp^{s}}^{(1+s)}
\\
\notag&=
X^{(pa)}
\{
\sum_{i=0}^{\min\{pc,pr\}}
\sum_{j=0}^{p^d-1}
\binom{j-pr-pc+2i}{i}
X^{(pr-i)}\mu_j^{(d)}Y^{(pc-i)}
\}
Y^{(pt)}\mu_{np^{b}-2pr+2pt}^{(1+b)}
\mu_{mp^{s}}^{(1+s)}
\\
\notag&=
\sum_{i=0}^{\min\{pc,pr\}}
\sum_{j=0}^{p^d-1}
\binom{j-pr-pc+2i}{i}
X^{(pa)}X^{(pr-i)}Y^{(pc-i)}
Y^{(pt)}\mu_{j+2(pc-i+pt)}^{(d)}\mu_{np^{b}-2pr+2pt}^{(1+b)}
\mu_{mp^{s}}^{(1+s)}
\\
\notag&=
\sum_{i=0}^{\min\{pc,pr\}}
\sum_{j=0}^{p^d-1}
\binom{j-pr-pc+2i}{i}
\binom{pa+pr-i}{pa}X^{(pa+pr-i)}\binom{pc+pt-i}{pt}
Y^{(pc-i+pt)}
\\
\notag&\hspace{3cm}
\mu_{j+2(pc-i+pt)}^{(d)}
\mu_{np^{b}-2pr+2pt}^{(1+b)}
\mu_{mp^{s}}^{(1+s)}
\\
\notag&=
\sum_{i=0}^{\min\{c,r\}}
\sum_{j=0}^{p^d-1}
\binom{j-pr-pc+2ip}{ip}
\binom{pa+pr-ip}{pa}X^{(pa+pr-ip)}\binom{pc+pt-ip}{pt}
Y^{(pc-ip+pt)}
\\
\notag&\hspace{3cm}
\mu_{j+2(pc-ip+pt)}^{(d)}
\mu_{np^{b}-2pr+2pt}^{(1+b)}
\mu_{mp^{s}}^{(1+s)}
\\
\notag&=
\sum_{i=0}^{\min\{c,r\}}
\sum_{j=0}^{p^{d-1}-1}
\binom{jp-pr-pc+2ip}{ip}
\binom{a+r-i}{a}X^{(p(a+r-i))}\binom{c+t-i}{t}
Y^{(p(c-i+t))}
\\
\notag&\hspace{3cm}
\mu_{jp+2p(c-i+t)}^{(d)}
\mu_{np^{b}-2pr+2pt}^{(1+b)}
\mu_{mp^{s}}^{(1+s)}
\\
\notag&=
\sum_{i=0}^{\min\{c,r\}}
\sum_{j=0}^{p^{d-1}-1}
\binom{j-a-c+2i}{i}
\binom{a+r-i}{a}
\binom{c+t-i}{t}
X^{(p(a+r-i))}Y^{(p(c-i+t))}
\\
\notag&\hspace{3cm}
\mu_{p(j+2(c-i+t))}^{(d)}
\mu_{np^{b}-2pr+2pt}^{(1+b)}
\mu_{mp^{s}}^{(1+s)}
\end{align}
avec
\begin{equation}
\begin{split}
&\mu_{p(j+2(c-i+t))}^{(d)}
 \mu_{np^{b}-2pr+2pt}^{(1+b)}
\mu_{mp^{s}}^{(1+s)}
\\
&=
\begin{cases}
\mu_{p(j+2(c-i+t))}^{(d)}
\quad\text{si  
$b=s=0$, ou bien  si $b>0$ et $s=0$
avec}
\\
\hspace{1cm}\text{
$j+2(c-i+t)\equiv
np^{b-1}-2r+2t\mod p^b$,
ou bien si  $b=0$ et $s>0$ avec}
\\
\hspace{1cm}\text{
$j+2(c-i+t)\equiv
np^{b-1}-2r+2t\mod p^s$
ou bien  finalement si $b, s>0$ avec \`a la fois }
\\
\hspace{1cm}\text{
$j+2(c-i+t)\equiv
np^{b-1}-2r+2t\mod p^b$
et}
\\
\hspace{1cm}\text{
$j+2(c-i+t)\equiv
np^{b-1}-2r+2t\mod p^s$}
\\
0
\quad\text{sinon}.
\end{cases}
\end{split}
\end{equation}

\subsection{}
\'Etendons maintenant lin\'eairement \`a $\Dist(G)$ l'endomorphisme
$\Fr'$ \eqref{defFrprime}  de $\Dist(T)$
par
\begin{equation}
X^{(a)}\mu_n^{(1+b)}Y^{(c)}
\mapsto
X^{(ap)}\Fr'(\mu_{n}^{(1+b)})Y^{(cp)}
\end{equation}
pour tous $a,b,c\in\bbN$ et $n\in\bbZ$. Alors

\begin{teo}\label{descrDist}
L'application
\begin{equation}\label{defaltphi} 
\phi : \Dist(G)\to
\sum_{a,b,c\in\bbN}
\sum_{n\in[0,p^{b}[}
\Bbb{F}_{p}
X^{(pa)}Y^{(pc)}\mu_{np^{b}}^{(1+b)}
\end{equation}
d\'efinie par
\begin{equation} 
\phi(\mu)  = \mu_0\Fr'(\mu)\mu_0=
\Fr'(\mu)\mu_0  \,\, pour\, tout\,\,  \mu \in  \Dist(G)
\end{equation}
est un isomorphisme de $\Bbb{F}_{p}$-alg\`ebres
tel  que
$\Dist(\Fr)\circ\phi=\id_{\Dist(G)}$.
\end{teo}
 
Comme le lecteur s'en assurera, l'utilisation de la lettre $\phi$ pour d\'esigner \eqref{defaltphi}  n'entra\^{i}ne pas de confusion. Reprenons  les notations de la preuve de  
 \ref{lemmu0Dist}
et supposons tout d'abord $b=s=0$.
On a
\begin{equation}
\begin{split}
&\phi(X^{(a)}Y^{(c)})
\phi(X^{(r)}Y^{(t)})
=
X^{(ap)}Y^{(cp)}\mu_0
X^{(rp)}Y^{(tp)}\mu_0
\\
=&
\sum_{i=0}^{\min\{c,r\}}
\sum_{j=0}^{p^{l-1}-1}
\binom{j-r-c+2i}{i}
\binom{a+r-i}{a}
\binom{c+t-i}{t}
X^{(p(a+r-i))}Y^{(p(c-i+t))}
\mu_{p(j+2(c-i+t))}^{(l)}
\end{split}
\end{equation}
alors que
\begin{equation}
\begin{split}
&\phi
(X^{(a)}Y^{(c)}X^{(r)}Y^{(t)})
\\
&=\phi(\sum_{i=0}^{\min\{c,r\}}
\sum_{j=0}^{p^{l-1}-1}
\binom{j-r-c+2i}{i}
\binom{a+r-i}{a}\binom{c+t-i}{t}
X^{(a+r-i)}
Y^{(c-i+t)}
\mu_{j+2(c-i+t)}^{(l-1)})
\\
&=
\sum_{i=0}^{\min\{c,r\}}
\sum_{j=0}^{p^{l-1}-1}
\binom{j-r-c+2i}{i}
\binom{a+r-i}{a}\binom{c+t-i}{t}
X^{p(a+r-i)}
Y^{p(c-i+t)}
\mu_0\Fr'(\mu_{j+2(c-i+t)}^{(l-1)})
\\
&=
\sum_{i=0}^{\min\{c,r\}}
\sum_{j=0}^{p^{l-1}-1}
\binom{j-a-c+2i}{i}
\binom{a+r-i}{a}
\binom{c+t-i}{t}
X^{p(a+r-i)}Y^{p(c-i+t)}
\mu_{p(j+2(c-i+t))}^{(l)}.
\end{split}
\end{equation}
Il s'ensuit bien que
\begin{equation}
\phi(X^{(a)}Y^{(c)})
\phi(X^{(r)}Y^{(t)})
=
\phi(X^{(a)}Y^{(c)}X^{(r)}Y^{(t)}).
\end{equation}

Si maintenant $b>0$ et $s=0$,
\begin{multline*}
\phi(X^{(a)}Y^{(c)}\mu_{nb^{b-1}}^{(b)})
\phi(X^{(r)}Y^{(t)})
=
X^{(ap)}Y^{(cp)}\mu_{np^b}^{(1+b)}
X^{(rp)}Y^{(tp)}\mu_0
\\
=
\begin{cases}
\sum_{i=0}^{\min\{c,r\}}
\sum_{j=0}^{p^{l-1}-1}
\binom{j-r-c+2i}{i}
\binom{a+r-i}{a}
\binom{c+t-i}{t}
X^{(p(a+r-i))}Y^{(p(c-i+t))}
\mu_{p(j+2(c-i+t))}^{(l)}
\\
\hspace{4cm}\text{si
$j+2(c-i+t)\equiv np^{b-1}-2r+2t\mod p^b$}
\\
0
\quad\text{sinon}
\end{cases}
\end{multline*}
alors que
\begin{equation}
\begin{split}
&\phi(X^{(a)}Y^{(c)}\mu_{np^{b-1}}^{(b)}X^{(r)}Y^{(t)})\\
&=
\phi(\sum_{i=0}^{\min\{c,r\}}
\sum_{j=0}^{p^{l-1}-1}
\binom{j-r-c+2i}{i}
\binom{a+r-i}{a}\binom{c+t-i}{t}
X^{(a+r-i)}
Y^{(c-i+t)}
\mu_{j+2(c-i+t)}^{(l-1)}\mu_{np^{b-1}-2r+2t}^{(b)})
\\
&=
\begin{cases}
\sum_{i=0}^{\min\{c,r\}}
\sum_{j=0}^{p^{l-1}-1}
\binom{j-r-c+2i}{i}
\binom{a+r-i}{a}\binom{c+t-i}{t}
X^{p(a+r-i)}
Y^{p(c-i+t)}
\mu_0\Fr'(\mu_{j+2(c-i+t)}^{(l-1)})
\\
\hspace{5cm}\text{si
$j+2(c-i+t)\equiv np^{b-1}-2r+2t\mod p^b$}
\\
0
\quad\text{sinon}
\end{cases}
\\
&=
\begin{cases}
\sum_{i=0}^{\min\{c,r\}}
\sum_{j=0}^{p^{l-1}-1}
\binom{j-a-c+2i}{i}
\binom{a+r-i}{a}
\binom{c+t-i}{t}
X^{p(a+r-i)}Y^{p(c-i+t)}
\mu_{p(j+2(c-i+t))}^{(l)}
\\
\hspace{5cm}\text{si
$j+2(c-i+t)\equiv np^{b-1}-2r+2t\mod p^b$}
\\
0
\quad\text{sinon},
\end{cases}
\end{split}
\end{equation}
et donc
\begin{equation}
\phi(X^{(a)}Y^{(c)}\mu_{nb^{b-1}}^{(b)})
\phi(X^{(r)}Y^{(t)})
=
\phi(X^{(a)}Y^{(c)}\mu_{nb^{b-1}}^{(b)}X^{(r)}Y^{(t)}).
\end{equation}

On traite de la m\^eme fa\c{c}on les cas restants.

\subsection{}
On laisse au lecteur le soin de v\'erifier qu'un raffinement des arguments donn\'es ci-dessus permet de voir d'une part que pour tout $r>0$, l'ensemble 
$\sum_{a,b,c\in[0,p^{r}[}
\sum_{n\in[0,p^{b}[}
\Bbb{F}_{p}
X^{(pa)}Y^{(pc)}\mu_{np^{b}}^{(1+b)}$
est une sous-$\Bbb{F}_{p}$-alg\`ebre de 
$\mu_0\Dist(G_{r+1})\mu_0$  et d'autre part que  l'application    $\phi : \Dist(G_{r})\to
\sum_{a,b,c\in[0,p^{r}[} \sum_{n\in[0,p^{b}[} \Bbb{F}_{p} X^{(pa)}Y^{(pc)}\mu_{np^{b}}^{(1+b)}$  analogue \`a   \eqref{defaltphi}  est un isomorphisme. Le lemme \ref{lemmu0Dist} et le th\'eor\`eme \ref{descrDist} se d\'eduisent alors de ces r\'esultats  par passage \`a la limite inductive suivant $r$.

  \section{Le cas  g\'en\'eral  }{}\label{ }
 
\subsection{}\label{}  
Les notations et hypoth\`eses g\'en\'erales   de l'introduction \ref{intro1} et \ref{intro3} sont d\'esormais en vigueur. On remarquera que celles sur $l$  implique que les ordres des    $q_i^2 \in  \cB$ sont tous \'egaux \`a $l$.
 Rappelons que  
$\lbr
K_i
\\
t
\rbr
=
\prod_{s=1}^t\frac{
K_iv_i^{-s+1}
-
K_i^{-1}
v_i^{s-1}
}
{v_i^s
-
v_i^{-s}
}
\in
U_{\bbQ(v)}$,
$i\in[1,\ell]$,
$t\in\bbN$,
appartient \`a
$U$.
Comme nous l'avons d\'ej\`a sous-entendu plus haut \eqref{defbinomK}, nous noterons 
$\lbr
K_i
\\
t
\rbr
\otimes1\in
U\otimes_\cA\cB$
simplement par le   symbole
$\lbr
K_i
\\
t
\rbr
$.

\subsection{}\label{} 
Posons, pour $i\in[1,\ell]$, $j\in[0,2l[$
\begin{equation}
\kappa_{ij}=
\frac{1}{2l}\sum_{r=0}^{2l-1}
q^{\frac{-jr}{2}}K_i^r\in
U_q\,\,\, ; \,\,\, \kappa=\prod_{i=1}^\ell\kappa_{i0}.
\end{equation} 
Lorsque $\ell=1$, ces \'el\'ements correspondent donc \`a ceux  consid\'er\'es dans 
\eqref{defkappa'j}.

\begin{lem}\label{lemgen} 
\begin{itemize}

\item[{\rm (i)}]
L'\'el\'ement $\kappa$ est l'unique idempotent  non nul de la 
sous-alg\`ebre
$u_\cB^0$
de
$U_\cB$ 
engendr\'ee par les 
$K_i$, $1\leq i\leq\ell$,
tel que
$K_i\kappa=\kappa$
pour tout $i\in[1,\ell]$.

\item[{\rm (ii)}]
 Pout tous $ i, j\in [1,\ell]$ et  $r\in\bbZ$, on a  
\begin{equation}\label{}
E_i^{(l r)}\kappa_{j0}=
\kappa_{j0} E_i^{(l r)}
\quad\text{et}\quad
F_i^{(l r)}\kappa_{j0}=
\kappa_{j0} F_i^{(l r)}.
\end{equation} 

\item[{\rm (iii)}]
Chaque $\kappa_{ij}$, avec
$i\in[1,\ell], j\in[1,2l[$, est un  idempotent,
et les  $\prod_{i=1}^\ell\kappa_{ij
}$,
$j\in[0,2l[$, 
forment  une d\'ecomposition de   $1$  en idempotents orthogonaux deux \`a deux.

\end{itemize}
\end{lem}

Les assertions  (i)  et  (ii)  d\'ecoulent du cas $\ell=1$  car les
$K_i$,
$i\in[1,\ell]$
commutent entre eux.

L'assertion   (iii)  d\'ecoule de   \ref{integ}


\subsection{}\label{}

Soient
$\bU_\cB^+=
\cB[X_i^{(r)}\mid
i\in I, r\in\bbN]$
et
$\bU_\cB^-=
\cB[Y_i^{(r)}\mid
i\in I, r\in\bbN]$ comme en \ref{intro4}.
Rappelons que si $H_i=[X_i,Y_i]$, chaque
$\binom{H_i}{r}=
\frac{H_i(H_i-1)\dots(H_i-r+1)}{r!}$, pour $i\in I$, $r\in\bbN$,  est un \'el\'ement de $\bU
\otimes_\bbZ\bbQ$,
 qui appartient en fait \`a  $\bfU$ et que l'on d\'esignera abusivement par le m\^eme symbole, de m\^eme que ses images canoniques dans   $\bU_\cA= \bU \otimes_\bbZ \cA$,  $\bU_\cB= \bU \otimes_\bbZ \cB$, $\cdots$.   
Si
 $\bU_\cB^0=
\cB[
\binom{H_i}{r}\mid
i\in I,
r\in\bbN]$,
on dispose d'un isomorphisme  $\cB$-lin\'eaire
$\bU_\cB^-\otimes_\cB
\bU_\cB^0\otimes_\cB
\bU_\cB^+\to
\bU_\cB$
donn\'e par la multiplication.
Soit
$\bU_\cB^{\geq0}$
(resp.
$\bU_\cB^{\leq0}$)
l'image de
$\bU_\cB^0\otimes_\cB
\bU_\cB^+$ (resp.
$\bU_\cB^-\otimes_\cB
\bU_\cB^0$) par cette application. 
C'est une sous-$\cB$-alg\`ebre de
$\bU_\cB$. 

\begin{lem}
Il existe des applications multiplicatives
$
\cB$-lin\'eaires
$\phi^{\geq0}:\bU_\cB^{\geq0}\to
U_\cB^{\geq0}$
et
$\phi^{\leq0}:\bU_\cB^{\leq0}\to
U_\cB^{\leq0}$ 
telles que pour tous $i\in I$ et $r\in\bbN$,
\begin{equation}\label{}
\phi^{\geq0}(X_i^{(r)})=
E_i^{(rl)}\kappa,
\quad
\phi^{\geq0}(\binom{H_i}{r})=
\lbr
K_i
\\
rl
\rbr
\kappa
=
\phi^{\leq0}(\binom{H_i}{r}),
\quad
\phi^{\leq0}(Y_i^{(r)})=
F_i^{(rl)}\kappa.
\end{equation} 
\end{lem}


En effet, il existe 
d'apr\`{e}s
\cite[Th.1.2]{kl}
des homomorphismes
de   $\cB$-alg\`ebres
${'\phi}^{\geq0}: \bU_\cB^{\geq0}\to
U_\cB^{\geq0}$
et
${'\phi}^{\leq0}:\bU_\cB^{\leq0}\to
U_\cB^{\leq0}$
modulo
$(K_i^l-1\mid
i\in I)$
tels que pour tous 
$i\in I$,
et
$r\in\bbN$, on ait
\begin{equation}\label{}
'\phi^{\geq0}(X_i^{(r)})=
E_i^{(rl)},
\quad
'\phi^{\geq0}(\binom{H_i}{r})=
\lbr
K_i
\\
rl
\rbr
=
{'\phi}^{\leq0}(\binom{H_i}{r}),
\quad
'\phi^{\leq0}(Y_i^{(r)})=
F_i^{(rl)}.
\end{equation}
Comme $K_i^l$ est un \'el\'ement central (\cite[Lem.4.4]{L89}) dans  $U_\cB$  
 et comme
$\kappa_{i0}(K_i^l-1)=0$
pour tout $i$ gr\^ace \`a \ref{lemgen} {\rm{(iii)}}, l'assertion en d\'ecoule.

\subsection{}\label{}

Avec ces notations, le m\^eme argument que celui utilis\'e dans \cite[Th.1.4]{gk} prouve la
\begin{prop}\label{} 
Il existe une application $\cB$-lin\'eaire  multiplicative
$\phi:\bU_\cB\to
U_\cB$
prolongeant
$\phi^{\geq0}$ 
et 
$\phi^{\leq0}$,
qui, de plus, scinde l'homomorphisme de Frobenius quantique (\ref{Fr}) 
$\Fr : U_\cB\to
\bU_\cB$.
\end{prop}
Le th\'eor\`eme (\ref{teointrod}) est ainsi d\'emontr\'e.

\subsection{}\label{}

L'application de Frobenius quantique (\ref{Fr}) $\Fr$  est en r\'ealit\'e un homomorphisme de
 $\cB$-alg\`ebres de  Hopf  \cite[1.1, 1.3]{lus} v\'erifiant quelques compatibilit\'es suppl\'ementaires. 
La comultiplication sur  $U$ v\'erifie
\begin{equation}
\Delta(E_i^{(n)})=
\sum_{j=0}^nv_i^{j(n-j)}E_i^{(n-j)}K_i^j
\otimes
E_i^{(j)},
\Delta(F_i^{(n)})=
\sum_{j=0}^nv_i^{-j(n-j)}F_i^{(j)}\otimes
K_i^{-j}F_i^{(n-j)},
 \Delta(K_i)=K_i\otimes
 K_i
\end{equation}

pour tous  $i\in[1,\ell]$ et $n\in\bbN$.
On a donc
$(\Fr\otimes\Fr)\circ\Delta(E_i^{(n)})
=
\sum_{j=0}^nX_i^{(\frac{n-j}{l})}
\otimes
X_i^{(\frac{j}{l})}$,
quantit\'e qui s'annule sauf si   $l\vert n$, auquel cas 
\begin{equation}
(\Fr\otimes\Fr)\circ\Delta(E_i^{(nl)})
=
\sum_{j=0}^nX_i^{(n-j)}
\otimes
X_i^{(j)}
=
\Delta\circ\Fr(E_i^{(nl)})
\end{equation}
pout tous $n\in\bbN$ et $i\in[1,\ell]$.
De m\^eme
$(\Fr\otimes\Fr)\circ\Delta(F_i^{(n)})=\Delta\circ\Fr(F_i^{(n)})$.
On a \'egalement
$(\Fr\otimes\Fr)\circ\Delta(K_i)=
(\Fr\otimes\Fr)(K_i\otimes K_i)
=
1\otimes1
=
\Delta\circ\Fr(K_i)$.

Comme
 $U=\cA[E_i^{(n)}, F_i^{(n)}, K_i^{\pm1}\mid
i\in[1,\ell],
n\in\bbN]$,
on en d\'eduit donc que 
\begin{equation}
(\Fr\otimes\Fr)\circ\Delta=
\Delta\circ\Fr.
\end{equation}
La co\"unit\'e
$\varepsilon$ sur
$U$ annule tous les 
$E_i^{(n)}$ et $F_i^{(n)}$,
$n>0$, et envoie $K_i$ to 1.
On a donc
\begin{equation}
\Fr\circ\varepsilon=\varepsilon\circ\Fr.
\end{equation}
Finalement, l'antipode $S$ sur $U$ est donn\'ee par
$S(E_i^{(n)})=(-1)v_i^{n(n-1)}K_i^{-n}E_i^{(n)}$,
$S(F_i^{(n)})=(-1)v_i^{-n(n-1)}
F_i^{(n)}K_i^{n}$,
et
$S(K_i)=K_i^{-1}$.
On a donc
\begin{equation}
(\Fr\otimes\Fr)\circ
S=
\Delta\circ
S.
\end{equation}

Pour le scindage $\phi$ \eqref{defscind} la situation est toute autre. On a
$
\varepsilon\circ\phi
=
\phi\circ\varepsilon
$,
et
$S(\kappa)=\kappa$
car
$K_i^{2l}=1$.
Comme $K_i\kappa=\kappa$
pour tout 
$ i\in[1,\ell]$
gr\^ace \`a
\ref{lemgen} {\it{(i)}},
on a donc aussi
\begin{equation}
S\circ\phi
=
\phi\circ
S.
\end{equation}
 L'application $\phi$ ne commute n\'eanmmoins pas avec les comultiplications  :
$(\phi\otimes\phi)\circ\Delta(1)=\kappa\otimes\kappa\ne
\Delta(\kappa)
=
\Delta\circ\phi(1)$.
Cependant,
\begin{equation}\label{compadelta}
\begin{split}
\Delta(\kappa)(\kappa\otimes\kappa)
&=
(\kappa\otimes\kappa)
\prod_{i=1}^\ell
\{
\frac{1}{2l}\sum_{j=0}^{2l-1}K_i^j
\otimes
K_i^j
\}
\quad
\text{en calculant apr\`es extension des scalaires \`a   $\bbQ[q]$}
\\
&=
\kappa\otimes\kappa
\quad\text{car $\kappa$ est une mesure  invariante}.
\end{split} 
\end{equation}

Rappelons aussi
\cite[1.1]{lus}, \cite[3.1.3]{L93}
qu'il existe une involution $\Omega$   
et une anti-involution
$\Psi$   
de
 $U$ d\'efinies par 
\begin{equation}
\Omega(E_{i})=F_{i}, \,\, \Omega(F_{i})=E_{i},\,\, \Omega(K_{i})=K_{i}^{-1}, \,\,\Omega(v)=v
\end{equation}
\begin{equation}
 \Psi(E_{i})=E_{i},\,\, \Psi(F_{i})=F_{i}, \,\,\Psi(K_{i})=K_{i}^{-1}, \,\,\Psi(v)=v 
\end{equation}
telles que l'anti-morphisme 
$\Psi\circ\Omega=\Omega\circ\Psi$
\'echange 
$E_i^{(n)}$ et $F_i^{(n)}$
et fixe  $K_i$ pour tout $i\in[1,\ell]$.
Abr\'egeant  
$(\Omega\circ\Psi)\otimes_\cA{\text{Id}}_{\cB}$
en
$\Omega\circ\Psi$,
on a  
\begin{equation}
\Omega\circ\Psi(\kappa)=\kappa.
\end{equation}
\\
Si maintenant 
$\tau$
d\'esigne l'anti-involution de   Chevalley de
$\bU_\bbQ$ 
\'echangeant chaque 
$X_i$ avec
 $Y_i$, 
 on a
\begin{equation}\label{omegapsicommut}
\Omega\circ\Psi\circ\phi=
\phi\circ\tau.
\end{equation}


\subsection{}\label{}

Nous renvoyons ici le lecteur \`a \cite[\S23]{L93} pour tout ce que nous utiliserons concernant la $\bbQ(v)$-alg\`ebre quantique \emph{modifi\'ee}   associ\'ee aux donn\'ees de \ref{intro1}. La $\cA$-forme de cette alg\`ebre quantique modifi\'ee  
peut s'\'ecrire
$\dot{U}=
\coprod_{\lambda\in\Lambda}U^+1_\lambda
U^-=
\coprod_{\lambda\in\Lambda}
U^-1_\lambda
U^+$
avec la structure de 
$U$-bimodule 
donn\'ee dans
\cite[23.1.3]{L93}.
On a en particulier
$K_i1_\lambda=v_i^{\langle\lambda,\alpha_i^\vee\rangle}1_\lambda$ pour tous $
i\in[1,\ell]$
et $\lambda\in\Lambda$.

Soit
\begin{equation}
\psi:U\to\dot{U} \,\, ;\,\, x\mapsto
x1_0
\end{equation}
l'application $\cA$-lin\'eaire donn\'ee par l'action \`a gauche. On a  $\psi(\kappa)=1_0$.

Posons $\bar\cA=\cA/(v-1)\simeq\bbZ$ et soit  
\begin{equation}
\eta:U\otimes_\cA\bar\cA\to
\bU
\end{equation}
le morphisme de passage au quotient $U\otimes_\cA\bar\cA
\to
(U\otimes_\cA\bar\cA)
/(K_i-1\mid i\in[1,\ell])\simeq\bU$. 

McGerty a construit dans \cite[Prop. 3.4]{Mc} un scindage du morphisme de Frobenius qui, dans notre situation, s'interpr\^ete comme une application
\begin{equation}
c :\dot{U}\otimes_\cA{\bar\cA}
 \to
\dot{U}\otimes_\cA(\cA/(\Phi_l)).
\end{equation}

Il r\'esulte alors imm\'ediatement des d\'efinitions qu'on a un diagramme commutatif
\begin{equation}
\xymatrix{
\dot{U}\otimes_\cA{\bar\cA}
\ar[rr]^-{c}
&&
\dot{U}\otimes_\cA(\cA/(\Phi_l))\ 
\ar@{^(->}[r]
&
\dot{U}\otimes_\cA\cB
\\
U\otimes_\cA\bar\cA
\ar[u]^{
\psi\otimes_\cA\bar\cA}
\ar[r]_-{
\eta}
&
\bU\ 
\ar@{^(->}[r]
&
\bU_\cB
\ar[r]_-{\phi}
&
U_\cB.
\ar[u]_-{\psi\otimes_\cA\cB}
}
\end{equation}

 

  \section{Contraction}{}\label{ }

\subsection{}\label{}
Comme dans le cas modulaire
\cite{gk}
on peut contracter tout 
$U_\cB$-module 
en utilisant le scindage $\phi$.
Comme $\kappa$ est un idempotent,
tout 
$U_\cB$-module
$M$
admet une d\'ecomposition
$M=\kappa
M\oplus(1-\kappa)M$.
Voyant $\kappa M$ comme un $\kappa U_\cB\kappa$-module, on peut d\'efinir, via l'homomorphisme de $\cB$-alg\`ebres \eqref{defscind}  
$\phi : \bU_\cB 
\rightarrow  \kappa U_\cB\kappa$, une structure de  
$\bU_\cB$-module sur $\kappa M$  
annulant ainsi
$(1-\kappa)M$.
On notera parfois  $M^\phi$ cette nouvelle structure sur 
$\kappa M$ et  on l'appellera la \emph{contraction par Frobenius} de
$M$.
On \'ecrira
$x\bullet m=\phi(x)m$
pour tous
$x\in
\bU_{\cB}$ et
 $m\in\kappa M$.

\subsection{}\label{}

Soit $\Lambda$ comme ci-dessus le r\'eseau des poids de 
$A$.
Pour tout
$\lambda\in\Lambda$
on d\'efinit un homomorphisme de $\cA$-alg\`ebres
\begin{equation}
\chi_\lambda:U^0\to
\cA
\end{equation}
par
\begin{equation}
 K_i\mapsto
v_i^{\langle\lambda,\alpha_i^\vee\rangle}\,\, ; \,\,  \lbr
K_i
\\
r
\rbr\mapsto
\lbr
\langle\lambda,\alpha_i^\vee\rangle
\\
r
\rbr_i
=
\prod_{s=1}^{r}\frac{v_i^{\langle\lambda,\alpha_i^\vee\rangle-s+1}-v_i^{-\langle\lambda,\alpha_i^\vee\rangle+s-1}}{v_i^s-v_i^{-s}}
\end{equation}
pour tout $i$ et $r\in\bbN$
\cite[1.1]{APW91}.
On notera  encore  simplement
$\chi_\lambda\otimes_\cA {\rm{Id}}_{\cB}$
par $\chi_\lambda : U^0_{\cB}\to
\cB$. 

Pour tout $\lambda\in\Lambda$, 
on d\'efinit de m\^eme un homomorphisme de 
$
\cB$-alg\`ebres
\begin{equation}
\bar\chi_{\lambda}:\bU_{\cB}^0\to
\cB
\end{equation}
par
\begin{equation}
\left(\begin{array}{c}H_i \\n\end{array}\right) \mapsto
\left(\begin{array}{c}\langle\lambda,\alpha_i^\vee\rangle \\n\end{array}\right)
\end{equation}
pour tous $i$ et $n\in {\Bbb{N}}$. Le lecteur remarquera que cette notation est compatible, en un sens \'evident, avec celle de \eqref{chi}.
On a   
\begin{equation}
\bar\chi_{\lambda}\circ\Fr|_{U_\cB^0}
=
 \chi_{l\lambda}.
\end{equation}

 \subsection{}\label{modintpoids}

Soient $\{\alpha_i
\mid
i\in[1,\ell]\}$ l'ensemble des racines simples correspondant aux  
$E_i$, et  
$\alpha_i^\vee$ les co-racines correspondantes.
Soit  
$\Lambda_1=\{\lambda\in\Lambda\mid
\langle\lambda,\alpha_i^\vee\rangle\in[0,l[\ \, {\text{pour tout}}\,  i
\}$.
Pour tout $\lambda\in\Lambda$,
nous \'ecrirons
$\lambda=\lambda^0+l\lambda^1$
avec
$\lambda^0\in\Lambda_1$
et
$\lambda^1\in\Lambda$.

Notons $\cC_\cB$ la cat\'egorie des $U_\cB$-modules int\'egrables de   type $\bf1$ \cite[1.6]{APW91}.
On dira qu'un objet $M$ de $\cC_\cB$ {\emph{se d\'ecompose suivant ses poids}} si  $M=\coprod_{\lambda\in\Lambda}M_\lambda$ avec
$M_\lambda=\{
m\in M\mid
xm=\chi_\lambda(x)m\ \, \text{pour tout}\,
x\in U_\cB^0\}$.

\begin{prop}\label{}
Pour tout $\lambda\in\Lambda$, on a
\begin{equation}
 \chi_\lambda\circ\phi|_{\bU_{
\cB}^0}
=
\begin{cases}
\bar\chi_{\lambda^1}
&\text{si $\lambda\in l\Lambda$},
\\
0
&\text{sinon}.
\end{cases}
\end{equation}
En particulier, pour tout 
$M\in\cC_\cB$ 
se d\'ecomposant suivant ses poids, on a
\begin{equation}
M^\phi
=
\coprod_{\lambda\in\Lambda}M_{l\lambda},
\end{equation}
avec
$\bU_{
\cB}^0$
agissant sur
$M_{l\lambda}$
par
$\bar\chi_{\lambda}$.

\end{prop}

En effet, on a  
($q_i$ d\'esignant l'image de  
$v_i$dans  $\cB$)
\begin{equation}
\begin{split}
 \chi_\lambda(\kappa)
&=
 \chi_\lambda(\prod_{i=1}^\ell
\kappa_{i0})
=
\prod_{i=1}^\ell\chi_\lambda(\kappa_{i0})
=
\prod_{i=1}^\ell
\{\frac{1}{2}\sum_{j=0}^{l-1}(-1)^j
\lbr
\langle\lambda,
\alpha_i^\vee\rangle
\\
j
\rbr_{i}
(q_i^j+q_i^{-j+\langle\lambda,\alpha_i^\vee\rangle}
)\}
\\
&=
\begin{cases}
1
&\text{si $l|\langle\lambda,\alpha_i^\vee\rangle$
\,{\text{pour tout}}\,$ i$}
\\
0
&\text{sinon}.
\end{cases}
\end{split}
\end{equation}

Il s'ensuit que pour tout $m\in\bbN$ 
\begin{equation}
\begin{split}
 \chi_\lambda\circ\phi
\binom{H_i}{m}
&=
\chi_\lambda
(\lbr K_i
\\
ml
\rbr)
\chi_\lambda(\kappa)
=
\begin{cases}
\lbr\langle\lambda,\alpha_i^\vee\rangle
\\
ml\rbr_i=
\lbr
l\langle\lambda^1,\alpha_i^\vee\rangle\\
ml\rbr_i=
\binom{\langle\lambda^1,\alpha_i^\vee\rangle}{m}
&\text{si $\lambda\in l\Lambda$
}
\\
0
&\text{sinon} .
\end{cases}
\end{split}
\end{equation}
Comme
$\bar\chi_{{\lambda^1}}
\binom{H_i}{m}
=
\binom{\langle\lambda^1,\alpha_i^\vee\rangle}{m}$, la proposition s'ensuit.

\subsection{}\label{}

Pour un   $U_q$-module de dimension finie $M$ 
d\'esignons par
$M^{\Omega\Psi}$
son  $\bbQ(q)$-dual  
$M^*$
equip\'e de la structure de
$U_q$-module d\'efinie par
$xf=f((\Omega\circ\Psi)(x)\,?)$, pour 
$x\in U_q$
et
$f\in M^*$. Posons enfin $\bU_q
= \bU_{\cB}\otimes_\cB\bbQ(q)$.

\begin{prop}
Pour tout  
$U_q$-module de dimension finie 
$M$, 
il existe un isomorphisme de  
$\bU_q
$-modules
\begin{equation}
(M^\phi)^{\tau}\simeq
(M^{\Omega\Psi})^\phi.
\end{equation}
\end{prop}
 
Consid\'erons  l'application de restriction  
$M^{\Omega\Psi}\to(\kappa M)^*$, 
qui est  $\bbQ(q)$-lin\'eaire et   surjective. 
Comme $\Omega\circ\Psi(\kappa)=\kappa$, 
on a
$\kappa f=f(\kappa\,?)$ pour tout 
$f\in M^{\Omega\Psi}$.
Comme $\kappa$ annule
$(1-\kappa)M$,
l'application de restriction induit la bijection voulue  
$(M^{\Omega\Psi})^\phi\simeq
(M^\phi)^\tau$
gr\^ace \`a  \eqref{omegapsicommut}.

\subsection{}\label{}
Si $M$ est un $\bU_\cB$-module, on   notera 
$M^\Fr$ le   $U_\cB$-module dont l'espace sous-jacent est celui de $M$ 
et la structure de  $U_\cB$-module est  celle donn\'ee compos\'ee avec $\Fr: U_\cB\to\bU_\cB$.
On a donc
$(M^{\Fr})^\phi=M$.

\begin{lem}\label{Frphi}
Soit $V$ un
$U_\cB$-module int\'egrable de type {\bf{1}}  
 et $M$ un $\bU_\cB$-module de type fini sur $\cB$
 admettant une d\'ecomposition suivant ses poids. On a des isomorphismes de $\bU_\cB$-modules
\begin{equation}\label{compatwist}
(V\otimes M^\Fr)^\phi\simeq
V^\phi\otimes M\simeq
(M^\Fr\otimes V)^\phi.
\end{equation}
\end{lem}

Soient 
$z\in V$ et $m\in M^\Fr$.
Pour tous  $i\in[1,\ell]$ et $r\in\bbN$, on a, dans  $(V\otimes M^\Fr)^\phi$,
\begin{equation}
\begin{split}
X_i^{(r)}\bullet
(\kappa z\otimes m)
&=
\phi(X_i^{(r)})
(\kappa z\otimes m)
=
\phi(X_i^{(r)})
(\kappa z\otimes\kappa m)
\quad\text{car
$\kappa m=m$}
\\
&=
\Delta(E_i^{(rl)})\Delta(\kappa)
(\kappa z\otimes\kappa m)
\\
&=
\Delta(E_i^{(rl)})(\kappa z\otimes\kappa m)
\quad\text{par \eqref{compadelta}}
\\
&=
\sum_{j=0}^{rl}
(q_i^{j(rl-j)}E_i^{(rl-j)}\otimes E_i^{(j)})(\kappa z\otimes m)
\quad\text{avec $q_i=q^{d_i}$}
\\
&=
\sum_{j=0}^{r}
(E_i^{((r-j)l)}\kappa z)\otimes(X_i^{(j)} m)
\quad\text{car $m\in M^\Fr$}
\end{split}
\end{equation}
alors que, regardant  $\kappa z\otimes
m$ dans  $V^\phi\otimes M$, on a 
\begin{equation}
\begin{split}
X_i^{(r)}\bullet
(\kappa z\otimes m)
&=
\Delta(X_i^{(r)})(\kappa z\otimes m)
=
\sum_{j=0}^{r}
(X_i^{(r-j)}\otimes
X_i^{(j)})
(\kappa z\otimes m)
\\
&=
\sum_{j=0}^{r}
(X_i^{((r-j)l)}\kappa z)\otimes(X_i^{(j)} m).
\end{split}
\end{equation}
De m\^eme pour  l'action  de
$Y_i^{(r)}$, si bien que le premier isomorphisme de \eqref{compatwist} s'ensuit. Le second isomorphisme se traite de mani\`ere analogue.

\subsection{}\label{}

Revenons \`a  la situation modulaire   et rappelons qu'alors qu'en caract\'eristique nulle  tous les  $\bU_\bbQ$-modules de dimension finie sont   semi-simples \cite[Th. II.8]{Jac}, tel n'est pas le cas pour les $G$-modules.

Soient $\Lambda^+$ l'ensemble des poids dominants,  $\nabla(\lambda)$ le $G$-module induit standard construit \`a partir de $\lambda\in\Lambda^+$ (ceux-ci sont d\'efinis sur   $\bbZ$ et fournissent les    $\bU_\bbQ$-modules simples par changement de base) et $L(\nu)$ le $G$-module simple de plus haut poids   $\nu \in \Lambda^+$. 
Soit $\Lambda_1$ comme dans  \ref{modintpoids} mais avec
$l$ remplac\'e par  $p$.

On dit qu'un $G$-module $M$ de dimension finie admet une {\emph{bonne filtration}} (resp. une {\emph{bonne $p$-filtration}}) si et seulement s'il admet une filtration   par des $G$-sous-modules  dont les gradu\'es associ\'es sont de la forme $\nabla(\lambda)$ avec $\lambda\in\Lambda^+$ (resp. $L(\nu)\otimes\nabla(\mu)^\Fr$ avec
$\nu\in\Lambda_1$, $\mu\in\Lambda^+$).
   
Soit $h$ le nombre de Coxeter de   $G$.

\begin{prop}\label{bonne}
Supposons $p\geq2(h-1)$ et soit $M$ un $G$-module $M$ de dimension finie. Toute   bonne $p$-filtration sur $M$ induit une bonne filtration sur  $M^\phi$.     
\end{prop}

Comme le foncteur $V \rightarrow V^\phi$ est exact, on peut supposer que
$M=L(\nu)\otimes\nabla(\lambda)^\Fr$ pour certains
$\nu\in\Lambda_1$ et $\lambda\in\Lambda^+$.
On a alors
$(L(\nu)\otimes\nabla(\lambda)^\Fr)^\phi\simeq
L(\nu)^\phi\otimes\nabla(\lambda)$ gr\^ace au lemme \ref{Frphi}.
Si $L(\nu)^\phi$ admet une bonne filtration avec des sous-quotients  
$\nabla(\eta)$, alors
$L(\nu)^\phi\otimes\nabla(\lambda)$ admettra une filtration avec des  
sous-quotients 
$\nabla(\eta)\otimes\nabla(\lambda)$, lesquels admettent 
  une bonne filtration gr\^ace \`a  
\cite[Thm 1.-1]{m}.

Encore gr\^ace au  lemme \ref{Frphi}, on peut m\^eme supposer  
$M=L(\nu)$ pour un certain  $\nu\in\Lambda_1$.
N\'eanmmoins, si $p\eta$ est un poids de  $L(\nu)$,
et si  $\alpha_0^\vee$ est la plus haute coracine de $G$
et 
$\rho$ la demi-somme  des racines positives, on a  
$p\langle\eta+\rho,\alpha_0^\vee\rangle=
\langle
p\eta,\alpha_0^\vee\rangle
+p(h-1)\leq\langle\nu,\alpha_0^\vee\rangle+p(h-1)\leq
\langle
(p-1)\rho,\alpha_0^\vee\rangle
+
p(h-1)=(2p-1)(h-1)$,
et donc
$\langle\eta+\rho,\alpha_0^\vee\rangle\leq
(h-1)(2-\frac{1}{p})
<2(h-1)\leq p$.

Il d\'ecoule alors du  linkage principle fort que 
 $L(\nu)^\phi$
est une somme directe de 
$\nabla(\eta)$ avec  des  $\eta\in\Lambda^+$.

\subsection{}\label{}

Supposant \'etablie la conjecture de Lusztig 
concernant les caract\`eres irr\'eductibles de $G$,   laquelle est un th\'eor\`eme pour $p$ suffisamment grand, Parshall et Scott \cite[Th. 5.1]{PS}
montrent que pour $p\geq2(h-1)$  tout  
$\nabla(\lambda)$ avec $\lambda\in\Lambda^+$, admet une bonne $p$-filtration.   On d\'eduit donc de la proposition  \ref{bonne} le

\begin{cor}
Supposons \'etablie la conjecture de Lusztig
et $p\geq2(h-1)$.
Soit  $M$ un $G$-module $M$ de dimension finie. Toute  bonne filtration  sur $M$ en induit une sur $M^\phi$.
\end{cor}


\subsection{}\label{}
 Terminons  par un   r\'esultat g\'eom\'etrique concernant  la th\'eorie modulaire    et r\'epondant \`a la question  laiss\'ee ouverte \`a la fin de [5]. 
 Notons ici $X= G/B$ la vari\'{e}t\'{e} des drapeaux de $G$, 
 $\cO_X$ son faisceau structural  et $F$ l'homomorphisme de Frobenius absolu sur $X$. Le fait que  $\cO_{X}$ soit un facteur direct de $F_*\cO_{X}$ est bien connu et est central dans l'\'{e}tude de la g\'{e}om\'{e}trie de $X$.  Nous d\'eterminons   un autre facteur direct de $F_*\cO_{X}$.
\\ 

\begin{teo}\label{}
{\it{Supposons $p\geq h$.
Le faisceau inversible  
$\cL_{X}(-\rho)$
associ\'e au 
 $B$-module $k$ de dimension 1   d\'efini par
$-\rho$
est un facteur direct de
$F_*\cO_X$.}}
\end{teo}

Posons  $\bar{X}=G/G_1B$ 
et  $q: X \to \bar{X}$ l'homomorphisme canonique. On va montrer, de mani\`ere \'equivalente  que le faisceau inversible
 $\cL_{\bar{X}}(-p\rho)$
sur $\bar{X}$ associ\'e au
$G_1B$-module de dimension 1 d\'efini par $-p\rho$ 
est un facteur direct de $q_*\cO_X$.
Plus pr\'ecis\'ement encore, soit
$L((p-2)\rho)$ le $G$-module simple de plus haut poids  
$(p-2)\rho$. On va alors montrer que  
$L((p-2)\rho)\otimes\cL_{\bar{X}}(-p\rho)$ est un facteur direct de
$q_*\cO_X$.

On a 
$q_*\cO_X\simeq
\cL_{\bar{X}}(\hat\nabla(k))$
avec $\hat\nabla(k)$ le $G_1B$-module induit \`a partir du 
 $B$-module trivial
$k$.
Comme
$\hat\nabla(k)$ a pour module de t\^ete simple $L((p-2)\rho)\otimes
p(-\rho)$
\cite[II.9.6]{J},
posons
$\pi:
\hat\nabla(k)\to
L((p-2)\rho)\otimes
p(-\rho)$
l'homomorphisme  
de passage au quotient.
Ainsi 
\begin{equation} 
   \cL_{\bar{X}}(\pi):q_*\cO_X\to
\cL_{\bar{X}}(L((p-2)\rho)\otimes
p(-\rho))\simeq
L((p-2)\rho)\otimes
\cL_{\bar{X}}(-p\rho)
\end{equation}
est un \'epimorphisme. On va montrer, en suivant la strat\'egie de  
\cite{K94}, qu'il est scind\'e.

D\'esignons par 
$\ind_{G_1B}^G$ le foncteur d'induction de $G_1B$ \`a $G$ de la cat\'egorie des  $G_1B$-modules vers celle des  $G$-modules.
On remarque tout d'abord qu'on a un diagramme commutatif
\begin{equation}
\tiny{
\xymatrix{
{\bf{Hom}}_{\cO_{\bar{X}}
}(\cL_{\bar{X}}(
L((p-2)\rho)\otimes
p(-\rho)),q_*\cO_X)
\ar[r]^-{\sim}
\ar[d]_-{
{\bf{Hom}}_{\cO_{\bar{X}}
}(\cL_{\bar{X}}(
L((p-2)\rho)\otimes
p(-\rho)),\cL_{\bar{X}}(\pi))
}
&
\ind_{G_1B}^G
(L((p-2)\rho)^*\otimes
p\rho\otimes\hat\nabla(k))
\ar[d]^-{
\ind_{G_1B}^G
(L((p-2)\rho)^*\otimes
p\rho\otimes\pi)}
\\
{\bf{Hom}}_{\cO_{\bar{X}}
}(\cL_{\bar{X}}(
L((p-2)\rho)\otimes
p(-\rho)),\cL_{\bar{X}}(
L((p-2)\rho)\otimes
p(-\rho)))
\ar[r]^-{\sim}
&
\ind_{G_1B}^G
(L((p-2)\rho)^*\otimes
p\rho\otimes
L((p-2)\rho)\otimes
p(-\rho)).
}
}
\end{equation}
D'autre part, on a des isomorphismes canoniques 
\begin{equation}
\ind_{G_1B}^G
(L((p-2)\rho)^*\otimes
p\rho\otimes\hat\nabla(k))
\simeq
L((p-2)\rho)^*\otimes
\ind_{G_1B}^G
(
\hat\nabla(p\rho))
\simeq
L((p-2)\rho)^*\otimes\nabla(p\rho)
\end{equation}
et
\begin{equation}
\ind_{G_1B}^G
(L((p-2)\rho)^*\otimes
p\rho\otimes
L((p-2)\rho)\otimes
p(-\rho))
\simeq
L((p-2)\rho)^*\otimes
L((p-2)\rho).
\end{equation}
Posons maintenant
$\pi'=\rho\otimes\pi:
\hat\nabla(p\rho)\to
L((p-2)\rho)$. On est alors ramen\'e \`a montrer que
 $\ind_{G_1B}^G(
 \pi')\ne0$ car
$L((p-2)\rho)$ est simple.
Pour ce faire, il suffit  de montrer que
$\ind_{G_1B}^G(\ker\pi') \subsetneqq \nabla(p\rho)$.
Mais le seul  facteur de composition de $\hat\nabla(p\rho)$,   comme 
$G_1B$-module, qui pourrait contribuer \`a produire  
 $L((p-2)\rho)$
une fois appliqu\'e 
$\rR^\bullet\ind_{G_1B}^G$
est de la forme 
$L((p-2)\rho)\otimes
p\mu$,
$\mu\in\Lambda$.
D'autre part
$
[\hat\nabla(p\rho):
L((p-2)\rho)\otimes
p\mu]\ne0$ 
si et seulement si 
$\mu=0$,
auquel cas 
$
[\hat\nabla(p\rho):
L((p-2)\rho)]=1$
\cite[II,9.16]{J}. 
Il s'ensuit que 
$\ind_{G_1B}^G(\ker\pi') \subsetneqq  \nabla(p\rho)$.
\\

\begin{cor}\label{}
 Supposons $p\geq h$.
 \begin{itemize}
\item[{\rm (i)}]
La multiplicit\'e de 
$L((p-2)\rho)$ dans  $\nabla(p\rho)$
est \'egale \`a  1  et celui-ci apparait comme module de t\^ete dans   
$\nabla(p\rho)$.
\item[{\rm (ii)}]
Pour tout $i>0$, on a
$\rH^i(\bar{X},\cL_{\bar{X}}(\ker\pi'))=0$. 
\item[{\rm (iii)}]
Pour tout $r>0$, $\cL_X(-\rho)$ est un facteur direct de $F^r_*\cO_{X}$.
\end{itemize}
\end{cor}
 

\end{document}